
\magnification=1200
\hsize=13.50cm
\vsize=18cm
\hoffset=-2mm
\voffset=.8cm
\parindent=12pt   \parskip=0pt
\hfuzz=1pt

\pretolerance=500 \tolerance=1000  \brokenpenalty=5000

\catcode`\@=11

\font\eightrm=cmr8
\font\eighti=cmmi8
\font\eightsy=cmsy8
\font\eightbf=cmbx8
\font\eighttt=cmtt8
\font\eightit=cmti8
\font\eightsl=cmsl8
\font\sevenrm=cmr7
\font\seveni=cmmi7
\font\sevensy=cmsy7
\font\sevenbf=cmbx7

\font\sixrm=cmr6
\font\sixi=cmmi6
\font\sixsy=cmsy6
\font\sixbf=cmbx6

\font\douzebf=cmbx10 at 12pt

\font\twelvebf=cmbx10 at 12pt

\font\tencal=eusm10

\font\sevencal=eusm7

\font\fivecal=eusm5
\newfam\calfam
\textfont\calfam=\tencal
\scriptfont\calfam=\sevencal
\scriptscriptfont\calfam=\fivecal
\def\cal#1{{\fam\calfam\relax#1}}

\skewchar\eighti='177 \skewchar\sixi='177
\skewchar\eightsy='60 \skewchar\sixsy='60

\def\tenpoint{%
   \textfont0=\tenrm \scriptfont0=\sevenrm
   \scriptscriptfont0=\fiverm
   \def\rm{\fam\z@\tenrm}%
   \textfont1=\teni  \scriptfont1=\seveni
   \scriptscriptfont1=\fivei
   \def\oldstyle{\fam\@ne\teni}\let\old=\oldstyle
   \textfont2=\tensy \scriptfont2=\sevensy
   \scriptscriptfont2=\fivesy
   \textfont\itfam=\tenit
   \def\it{\fam\itfam\tenit}%
   \textfont\slfam=\tensl
   \def\sl{\fam\slfam\tensl}%
   \textfont\bffam=\tenbf
   \scriptfont\bffam=\sevenbf
   \scriptscriptfont\bffam=\fivebf
   \def\bf{\fam\bffam\tenbf}%
   \textfont\ttfam=\tentt
   \def\tt{\fam\ttfam\tentt}%
   \abovedisplayskip=12pt plus 3pt minus 9pt
   \belowdisplayskip=\abovedisplayskip
   \abovedisplayshortskip=0pt plus 3pt
   \belowdisplayshortskip=4pt plus 3pt
   \smallskipamount=3pt plus 1pt minus 1pt
   \medskipamount=6pt plus 2pt minus 2pt
   \bigskipamount=12pt plus 4pt minus 4pt
   \normalbaselineskip=12pt
   \setbox\strutbox=\hbox{\vrule height8.5pt depth3.5pt width0pt}%
   \let\bigf@nt=\tenrm
   \let\smallf@nt=\sevenrm
   \normalbaselines\rm}

\def\eightpoint{%
   \textfont0=\eightrm \scriptfont0=\sixrm
   \scriptscriptfont0=\fiverm
   \def\rm{\fam\z@\eightrm}%
   \textfont1=\eighti  \scriptfont1=\sixi
   \scriptscriptfont1=\fivei
   \def\oldstyle{\fam\@ne\eighti}\let\old=\oldstyle
   \textfont2=\eightsy \scriptfont2=\sixsy
   \scriptscriptfont2=\fivesy
   \textfont\itfam=\eightit
   \def\it{\fam\itfam\eightit}%
   \textfont\slfam=\eightsl
   \def\sl{\fam\slfam\eightsl}%
   \textfont\bffam=\eightbf
   \scriptfont\bffam=\sixbf
   \scriptscriptfont\bffam=\fivebf
   \def\bf{\fam\bffam\eightbf}%
   \textfont\ttfam=\eighttt
   \def\tt{\fam\ttfam\eighttt}%
   \abovedisplayskip=9pt plus 3pt minus 9pt
   \belowdisplayskip=\abovedisplayskip
   \abovedisplayshortskip=0pt plus 3pt
   \belowdisplayshortskip=3pt plus 3pt
   \smallskipamount=2pt plus 1pt minus 1pt
   \medskipamount=4pt plus 2pt minus 1pt
   \bigskipamount=9pt plus 3pt minus 3pt
   \normalbaselineskip=9pt
   \setbox\strutbox=\hbox{\vrule height7pt depth2pt width0pt}%
   \let\bigf@nt=\eightrm
   \let\smallf@nt=\sixrm
   \normalbaselines\rm}
\tenpoint


\def\pc#1{\bigf@nt#1\smallf@nt}

\catcode`\;=\active
\def;{\relax\ifhmode\ifdim\lastskip>\z@\unskip\fi \kern\fontdimen2
 -1.2 \fontdimen3 \string;}

\catcode`\:=\active
\def:{\relax\ifhmode\ifdim\lastskip>\z@\unskip\fi\penalty\@M\
\fi\string:}

\catcode`\!=\active
\def!{\relax\ifhmode\ifdim\lastskip>\z@ \unskip\fi\kern\fontdimen2
 -1.1 \fontdimen3 \string!}

\catcode`\?=\active
\def?{\relax\ifhmode\ifdim\lastskip>\z@ \unskip\fi\kern\fontdimen2
 -1.1 \fontdimen3 \string?}

\frenchspacing

\newtoks\auteurcourant
\auteurcourant={\hfil}

\newtoks\titrecourant
\titrecourant={\hfil}

\newtoks\hautpagetitre
\hautpagetitre={\hfil}

\newtoks\baspagetitre
\baspagetitre={\hfil}

\newtoks\hautpagegauche
\hautpagegauche={\eightpoint\rlap{\folio}\hfil\the
\auteurcourant\hfil}

\newtoks\hautpagedroite
\hautpagedroite={\eightpoint\hfil\the
\titrecourant\hfil\llap{\folio}}

\newtoks\baspagegauche
\baspagegauche={\hfil}

\newtoks\baspagedroite
\baspagedroite={\hfil}

\newif\ifpagetitre
\pagetitretrue


\headline={\ifpagetitre\the\hautpagetitre
\else\ifodd\pageno\the\hautpagedroite\else\the
\hautpagegauche\fi\fi}

\footline={\ifpagetitre\the\baspagetitre\else
\ifodd\pageno\the\baspagedroite\else\the
\baspagegauche\fi\fi
\global\pagetitrefalse}

\def\raggedbottom{\topskip 10pt plus 36pt\r@ggedbottomtrue}

\def\point{\raise.2ex\hbox{\douzebf .}}
\def\pointir{\unskip . --- \ignorespaces}


\def\Medbreak{\vskip-\lastskip\medbreak}

\def\rem#1\endrem{%
\Medbreak {\it#1\unskip} : }


\long\def\th#1 #2\enonce#3\endth{%
\Medbreak {\pc#1} {#2\unskip}\pointir{\it #3}\medskip}


\long\def\thm#1 #2\enonce#3\endthm{%
\Medbreak {\pc#1} {#2\unskip}\pointir{\it #3}\medskip}

\def\decale#1{\smallbreak\hskip 28pt\llap{#1}\kern 5pt}

\def\decaledecale#1{\smallbreak\hskip 34pt\llap{#1}\kern 5pt}

\let\@ldmessage=\message

\def\message#1{{\def\pc{\string\pc\space}%
\def\'{\string'}\def\`{\string`}%
\def\^{\string^}\def\"{\string"}%
\@ldmessage{#1}}}

\def\({{\rm (}}

\def\){{\rm )}}


\def\up#1{\raise 1ex\hbox{\smallf@nt#1}}

\def\diagram#1{\def\normalbaselines{\baselineskip=0pt
\lineskip=5pt}\matrix{#1}}


\def\longmaprightover#1#2{\smash{\mathop{\hbox to#2
{\rightarrowfill}}\limits^{\scriptstyle#1}}}

\def\longmapleftover#1#2{\smash{\mathop{\hbox to#2
{\leftarrowfill}}\limits^{\scriptstyle#1}}}

\def\longmaprightunder#1#2{\smash{\mathop{\hbox to#2
{\rightarrowfill}}\limits_{\scriptstyle#1}}}

\def\longmapleftunder#1#2{\smash{\mathop{\hbox to#2
{\leftarrowfill}}\limits_{\scriptstyle#1}}}

\def\longhookrightarrowover#1#2{\smash{\mathop{\lhook\joinrel
\mathrel{\hbox to #2{\rightarrowfill}}}\limits^{\scriptstyle#1}}}

\def\longhookrightarrowunder#1#2{\smash{\mathop{\lhook\joinrel
\mathrel{\hbox to #2{\rightarrowfill}}}\limits_{\scriptstyle#1}}}

\def\longhookleftarrowover#1#2{\smash{\mathop{{\hbox to #2
{\leftarrowfill}}\joinrel\kern -0.9mm\mathrel\rhook}
\limits^{\scriptstyle#1}}}

\def\longhookleftarrowunder#1#2{\smash{\mathop{{\hbox to #2
{\leftarrowfill}}\joinrel\kern -0.9mm\mathrel\rhook}
\limits_{\scriptstyle#1}}}

\def\longtwoheadrightarrowover#1#2{\smash{\mathop{{\hbox to #2
{\rightarrowfill}}\kern -3.25mm\joinrel\mathrel\rightarrow}
\limits^{\scriptstyle#1}}}

\def\longtwoheadrightarrowunder#1#2{\smash{\mathop{{\hbox to #2
{\rightarrowfill}}\kern -3.25mm\joinrel\mathrel\rightarrow}
\limits_{\scriptstyle#1}}}

\def\longtwoheadleftarrowover#1#2{\smash{\mathop{\joinrel\mathrel
\leftarrow\kern -3.8mm{\hbox to #2{\leftarrowfill}}}
\limits^{\scriptstyle#1}}}

\def\longtwoheadleftarrowunder#1#2{\smash{\mathop{\joinrel\mathrel
\leftarrow\kern -3.8mm{\hbox to #2{\leftarrowfill}}}
\limits_{\scriptstyle#1}}}


\def\longmapsto#1{\mapstochar\mathrel{\joinrel \kern-0.2mm\hbox to
#1mm{\rightarrowfill}}}

\def\og{\leavevmode\raise.3ex\hbox{$\scriptscriptstyle
\langle\!\langle\,$}}
\def\fg{\leavevmode\raise.3ex\hbox{$\scriptscriptstyle
\,\rangle\!\rangle$}}

\def\section#1#2{\vskip 5mm {\bf {#1}. {#2}}\vskip 5mm}
\def\subsection#1#2{\vskip 3mm {\it #2}\vskip 3mm}

\catcode`\@=12

\showboxbreadth=-1  \showboxdepth=-1


\message{`lline' & `vector' macros from LaTeX}

\def\Grille{\setbox13=\vbox to 5\unitlength{\hrule width 109mm \vfill}
\setbox13=\vbox to 65mm
{\offinterlineskip\leaders\copy13\vfill\kern-1pt\hrule}
\setbox14=\hbox to 5\unitlength{\vrule height 65mm\hfill}
\setbox14=\hbox to 109mm{\leaders\copy14\hfill\kern-2mm \vrule height
65mm}
\ht14=0pt\dp14=0pt\wd14=0pt \setbox13=\vbox to 0pt
{\vss\box13\offinterlineskip\box14} \wd13=0pt\box13}

\def\rule(#1,#2)\dir(#3,#4)\long#5{%
\noalign{\leftput(#1,#2){\lline(#3,#4){#5}}}}

\def\arrow(#1,#2)\dir(#3,#4)\length#5{%
\noalign{\leftput(#1,#2){\vector(#3,#4){#5}}}}

\def\put(#1,#2)#3{\noalign{\setbox1=\hbox{%
\kern #1\unitlength \raise #2\unitlength\hbox{$#3$}}%
\ht1=0pt \wd1=0pt \dp1=0pt\box1}}

\catcode`@=11

\def\{{\relax\ifmmode\lbrace\else$\lbrace$\fi}

\def\}{\relax\ifmmode\rbrace\else$\rbrace$\fi}

\def\newcount{\alloc@0\count\countdef\insc@unt}

\def\newdimen{\alloc@1\dimen\dimendef\insc@unt}

\def\newwrite{\alloc@7\write\chardef\sixt@@n}

\newwrite\@unused
\newcount\@tempcnta
\newcount\@tempcntb
\newdimen\@tempdima
\newdimen\@tempdimb
\newbox\@tempboxa

\def\@spaces{\space\space\space\space}

\def\@whilenoop#1{}

\def\@whiledim#1\do #2{\ifdim #1\relax#2\@iwhiledim{#1\relax#2}\fi}

\def\@iwhiledim#1{\ifdim #1\let\@nextwhile=\@iwhiledim
\else\let\@nextwhile=\@whilenoop\fi\@nextwhile{#1}}

\def\@badlinearg{\@latexerr{Bad \string\line\space or \string\vector
\space argument}}

\def\@latexerr#1#2{\begingroup
\edef\@tempc{#2}\expandafter\errhelp\expandafter{\@tempc}%

\def\@eha{Your command was ignored.^^JType \space I <command> <return>
\space to replace it with another command,^^Jor \space <return> \space
to continue without it.}

\def\@ehb{You've lost some text. \space \@ehc}

\def\@ehc{Try typing \space <return> \space to proceed.^^JIf that
doesn't work, type \space X <return> \space to quit.}

\def\@ehd{You're in trouble here.  \space\@ehc}

\typeout{LaTeX error.  \space See LaTeX manual for explanation.^^J
\space\@spaces\@spaces\@spaces Type \space H <return> \space for
immediate help.}\errmessage{#1}\endgroup}

\def\typeout#1{{\let\protect\string\immediate\write\@unused{#1}}}

\font\tenln = line10
\font\tenlnw = linew10

\newdimen\@wholewidth
\newdimen\@halfwidth
\newdimen\unitlength

\unitlength =1pt

\def\thinlines{\let\@linefnt\tenln \let\@circlefnt\tencirc
\@wholewidth\fontdimen8\tenln \@halfwidth .5\@wholewidth}

\def\thicklines{\let\@linefnt\tenlnw \let\@circlefnt\tencircw
\@wholewidth\fontdimen8\tenlnw \@halfwidth .5\@wholewidth}

\def\linethickness#1{\@wholewidth #1\relax \@halfwidth .5
\@wholewidth}

\newif\if@negarg

\def\lline(#1,#2)#3{\@xarg #1\relax \@yarg #2\relax
\@linelen=#3\unitlength \ifnum\@xarg =0 \@vline \else \ifnum\@yarg =0
\@hline \else \@sline\fi \fi}

\def\@sline{\ifnum\@xarg< 0 \@negargtrue \@xarg -\@xarg \@yyarg
-\@yarg \else \@negargfalse \@yyarg \@yarg \fi
\ifnum \@yyarg >0 \@tempcnta\@yyarg \else \@tempcnta - \@yyarg \fi
\ifnum\@tempcnta>6 \@badlinearg\@tempcnta0 \fi
\setbox\@linechar\hbox{\@linefnt\@getlinechar(\@xarg,\@yyarg)}%
\ifnum \@yarg >0 \let\@upordown\raise \@clnht\z@
\else\let\@upordown\lower \@clnht \ht\@linechar\fi
\@clnwd=\wd\@linechar
\if@negarg \hskip -\wd\@linechar \def\@tempa{\hskip -2\wd \@linechar}
\else \let\@tempa\relax \fi
\@whiledim \@clnwd <\@linelen \do {\@upordown\@clnht\copy\@linechar
\@tempa \advance\@clnht \ht\@linechar \advance\@clnwd \wd\@linechar}%
\advance\@clnht -\ht\@linechar \advance\@clnwd -\wd\@linechar
\@tempdima\@linelen\advance\@tempdima -\@clnwd
\@tempdimb\@tempdima\advance\@tempdimb -\wd\@linechar
\if@negarg \hskip -\@tempdimb \else \hskip \@tempdimb \fi
\multiply\@tempdima \@m\@tempcnta \@tempdima \@tempdima \wd\@linechar
\divide\@tempcnta \@tempdima \@tempdima \ht\@linechar
\multiply\@tempdima \@tempcnta \divide\@tempdima \@m \advance\@clnht
\@tempdima
\ifdim \@linelen <\wd\@linechar \hskip \wd\@linechar
\else\@upordown\@clnht\copy\@linechar\fi}

\def\@hline{\ifnum \@xarg <0 \hskip -\@linelen \fi
\vrule height \@halfwidth depth \@halfwidth width \@linelen
\ifnum \@xarg <0 \hskip -\@linelen \fi}

\def\@getlinechar(#1,#2){\@tempcnta#1\relax
\multiply\@tempcnta 8\advance\@tempcnta -9
\ifnum #2>0 \advance\@tempcnta #2\relax
  \else\advance\@tempcnta -#2\relax\advance\@tempcnta 64 \fi
\char\@tempcnta}

\def\vector(#1,#2)#3{\@xarg #1\relax \@yarg #2\relax
\@linelen=#3\unitlength
\ifnum\@xarg =0 \@vvector \else \ifnum\@yarg =0 \@hvector \else
\@svector\fi \fi}

\def\@hvector{\@hline\hbox to 0pt{\@linefnt \ifnum \@xarg <0
\@getlarrow(1,0)\hss\else \hss\@getrarrow(1,0)\fi}}

\def\@vvector{\ifnum \@yarg <0 \@downvector \else \@upvector \fi}

\def\@svector{\@sline\@tempcnta\@yarg \ifnum\@tempcnta <0
\@tempcnta=-\@tempcnta\fi \ifnum\@tempcnta <5 \hskip -\wd\@linechar
\@upordown\@clnht \hbox{\@linefnt \if@negarg
\@getlarrow(\@xarg,\@yyarg) \else \@getrarrow(\@xarg,\@yyarg)
\fi}\else\@badlinearg\fi}

\def\@getlarrow(#1,#2){\ifnum #2 =\z@ \@tempcnta='33\else
\@tempcnta=#1\relax\multiply\@tempcnta \sixt@@n \advance\@tempcnta -9
\@tempcntb=#2\relax \multiply\@tempcntb \tw@ \ifnum \@tempcntb >0
\advance\@tempcnta \@tempcntb\relax \else\advance\@tempcnta
-\@tempcntb\advance\@tempcnta 64 \fi\fi \char\@tempcnta}

\def\@getrarrow(#1,#2){\@tempcntb=#2\relax \ifnum\@tempcntb < 0
\@tempcntb=-\@tempcntb\relax\fi \ifcase \@tempcntb\relax
\@tempcnta='55 \or \ifnum #1<3 \@tempcnta=#1\relax\multiply\@tempcnta
24 \advance\@tempcnta -6 \else \ifnum #1=3 \@tempcnta=49
\else\@tempcnta=58 \fi\fi\or \ifnum #1<3
\@tempcnta=#1\relax\multiply\@tempcnta 24 \advance\@tempcnta -3 \else
\@tempcnta=51\fi\or \@tempcnta=#1\relax\multiply\@tempcnta \sixt@@n
\advance\@tempcnta -\tw@ \else \@tempcnta=#1\relax\multiply\@tempcnta
\sixt@@n \advance\@tempcnta 7 \fi \ifnum #2<0 \advance\@tempcnta 64
\fi \char\@tempcnta}

\def\@vline{\ifnum \@yarg <0 \@downline \else \@upline\fi}

\def\@upline{\hbox to \z@{\hskip -\@halfwidth \vrule width
\@wholewidth height \@linelen depth \z@\hss}}

\def\@downline{\hbox to \z@{\hskip -\@halfwidth \vrule width
\@wholewidth height \z@ depth \@linelen \hss}}

\def\@upvector{\@upline\setbox\@tempboxa
\hbox{\@linefnt\char'66}\raise \@linelen \hbox to\z@{\lower
\ht\@tempboxa \box\@tempboxa\hss}}

\def\@downvector{\@downline\lower \@linelen \hbox to
\z@{\@linefnt\char'77\hss}}

\thinlines

\newcount\@xarg
\newcount\@yarg
\newcount\@yyarg
\newcount\@multicnt
\newdimen\@xdim
\newdimen\@ydim
\newbox\@linechar
\newdimen\@linelen
\newdimen\@clnwd
\newdimen\@clnht
\newdimen\@dashdim
\newbox\@dashbox
\newcount\@dashcnt
\catcode`@=12

\newbox\tbox
\newbox\tboxa

\def\leftzer#1{\setbox\tbox=\hbox to 0pt{#1\hss}%
\ht\tbox=0pt \dp\tbox=0pt \box\tbox}

\def\rightzer#1{\setbox\tbox=\hbox to 0pt{\hss #1}%
\ht\tbox=0pt \dp\tbox=0pt \box\tbox}

\def\centerzer#1{\setbox\tbox=\hbox to 0pt{\hss #1\hss}%
\ht\tbox=0pt \dp\tbox=0pt \box\tbox}

\def\leftput(#1,#2)#3{\setbox\tboxa=\hbox{%
\kern #1\unitlength \raise #2\unitlength\hbox{\leftzer{#3}}}%
\ht\tboxa=0pt \wd\tboxa=0pt \dp\tboxa=0pt\box\tboxa}

\def\rightput(#1,#2)#3{\setbox\tboxa=\hbox{%
\kern #1\unitlength \raise #2\unitlength\hbox{\rightzer{#3}}}%
\ht\tboxa=0pt \wd\tboxa=0pt \dp\tboxa=0pt\box\tboxa}

\def\centerput(#1,#2)#3{\setbox\tboxa=\hbox{%
\kern #1\unitlength \raise #2\unitlength\hbox{\centerzer{#3}}}%
\ht\tboxa=0pt \wd\tboxa=0pt \dp\tboxa=0pt\box\tboxa}

\unitlength=1mm

\expandafter\ifx\csname amssym.def\endcsname\relax \else
\endinput\fi
%
\expandafter\edef\csname amssym.def\endcsname{%
        \catcode`\noexpand\@=\the\catcode`\@\space}
\catcode`\@=11
%

\def\undefine#1{\let#1\undefined}
\def\newsymbol#1#2#3#4#5{\let\next@\relax
  \ifnum#2=\@ne\let\next@\msafam@\else
  \ifnum#2=\tw@\let\next@\msbfam@\fi\fi
  \mathchardef#1="#3\next@#4#5}
\def\mathhexbox@#1#2#3{\relax
  \ifmmode\mathpalette{}{\m@th\mathchar"#1#2#3}%
  \else\leavevmode\hbox{$\m@th\mathchar"#1#2#3$}\fi}
\def\hexnumber@#1{\ifcase#1 0\or 1\or 2\or 3
\or 4\or 5\or 6\or 7\or 8\or
  9\or A\or B\or C\or D\or E\or F\fi}

\font\tenmsa=msam10
\font\sevenmsa=msam7
\font\fivemsa=msam5
\newfam\msafam
\textfont\msafam=\tenmsa
\scriptfont\msafam=\sevenmsa
\scriptscriptfont\msafam=\fivemsa
\edef\msafam@{\hexnumber@\msafam}
\mathchardef\dabar@"0\msafam@39
\def\dashrightarrow{\mathrel{\dabar@\dabar@\mathchar"0
\msafam@4B}}
\def\dashleftarrow{\mathrel{\mathchar"0\msafam@4C
\dabar@\dabar@}}

\def\ulcorner{\delimiter"4\msafam@70\msafam@70 }
\def\urcorner{\delimiter"5\msafam@71\msafam@71 }
\def\llcorner{\delimiter"4\msafam@78\msafam@78 }
\def\lrcorner{\delimiter"5\msafam@79\msafam@79 }
\def\yen{{\mathhexbox@\msafam@55}}
\def\checkmark{{\mathhexbox@\msafam@58}}
\def\circledR{{\mathhexbox@\msafam@72}}
\def\maltese{{\mathhexbox@\msafam@7A}}

\font\tenmsb=msbm10
\font\sevenmsb=msbm7
\font\fivemsb=msbm5
\newfam\msbfam
\textfont\msbfam=\tenmsb
\scriptfont\msbfam=\sevenmsb
\scriptscriptfont\msbfam=\fivemsb
\edef\msbfam@{\hexnumber@\msbfam}
\def\Bbb#1{{\fam\msbfam\relax#1}}
\def\widehat#1{\setbox\z@\hbox{$\m@th#1$}%
  \ifdim\wd\z@>\tw@ em\mathaccent"0\msbfam@5B{#1}%
  \else\mathaccent"0362{#1}\fi}
\def\widetilde#1{\setbox\z@\hbox{$\m@th#1$}%
  \ifdim\wd\z@>\tw@ em\mathaccent"0\msbfam@5D{#1}%
  \else\mathaccent"0365{#1}\fi}
\font\teneufm=eufm10
\font\seveneufm=eufm7
\font\fiveeufm=eufm5
\newfam\eufmfam
\textfont\eufmfam=\teneufm
\scriptfont\eufmfam=\seveneufm
\scriptscriptfont\eufmfam=\fiveeufm
\def\frak#1{{\fam\eufmfam\relax#1}}

\csname amssym.def\endcsname

\expandafter\ifx\csname pre amssym.tex at\endcsname\relax \else
\endinput\fi
\expandafter\chardef\csname pre amssym.tex at\endcsname=\the
\catcode`\@
\catcode`\@=11
\begingroup\ifx\undefined\newsymbol \else\def\input#1
{\endgroup}\fi
\input amssym.def \relax
\newsymbol\boxdot 1200
\newsymbol\boxplus 1201
\newsymbol\boxtimes 1202
\newsymbol\square 1003
\newsymbol\blacksquare 1004
\newsymbol\centerdot 1205
\newsymbol\lozenge 1006
\newsymbol\blacklozenge 1007
\newsymbol\circlearrowright 1308
\newsymbol\circlearrowleft 1309
\undefine\rightleftharpoons
\newsymbol\rightleftharpoons 130A
\newsymbol\leftrightharpoons 130B
\newsymbol\boxminus 120C
\newsymbol\Vdash 130D
\newsymbol\Vvdash 130E
\newsymbol\vDash 130F
\newsymbol\twoheadrightarrow 1310
\newsymbol\twoheadleftarrow 1311
\newsymbol\leftleftarrows 1312
\newsymbol\rightrightarrows 1313
\newsymbol\upuparrows 1314
\newsymbol\downdownarrows 1315
\newsymbol\upharpoonright 1316
  
\newsymbol\downharpoonright 1317
\newsymbol\upharpoonleft 1318
\newsymbol\downharpoonleft 1319
\newsymbol\rightarrowtail 131A
\newsymbol\leftarrowtail 131B
\newsymbol\leftrightarrows 131C
\newsymbol\rightleftarrows 131D
\newsymbol\Lsh 131E
\newsymbol\Rsh 131F
\newsymbol\rightsquigarrow 1320
\newsymbol\leftrightsquigarrow 1321
\newsymbol\looparrowleft 1322
\newsymbol\looparrowright 1323
\newsymbol\circeq 1324
\newsymbol\succsim 1325
\newsymbol\gtrsim 1326
\newsymbol\gtrapprox 1327
\newsymbol\multimap 1328
\newsymbol\therefore 1329
\newsymbol\because 132A
\newsymbol\doteqdot 132B
  
\newsymbol\triangleq 132C
\newsymbol\precsim 132D
\newsymbol\lesssim 132E
\newsymbol\lessapprox 132F
\newsymbol\eqslantless 1330
\newsymbol\eqslantgtr 1331
\newsymbol\curlyeqprec 1332
\newsymbol\curlyeqsucc 1333
\newsymbol\preccurlyeq 1334
\newsymbol\leqq 1335
\newsymbol\leqslant 1336
\newsymbol\lessgtr 1337
\newsymbol\backprime 1038
\newsymbol\risingdotseq 133A
\newsymbol\fallingdotseq 133B
\newsymbol\succcurlyeq 133C
\newsymbol\geqq 133D
\newsymbol\geqslant 133E
\newsymbol\gtrless 133F
\newsymbol\sqsubset 1340
\newsymbol\sqsupset 1341
\newsymbol\vartriangleright 1342
\newsymbol\vartriangleleft 1343
\newsymbol\trianglerighteq 1344
\newsymbol\trianglelefteq 1345
\newsymbol\bigstar 1046
\newsymbol\between 1347
\newsymbol\blacktriangledown 1048
\newsymbol\blacktriangleright 1349
\newsymbol\blacktriangleleft 134A
\newsymbol\vartriangle 134D
\newsymbol\blacktriangle 104E
\newsymbol\triangledown 104F
\newsymbol\eqcirc 1350
\newsymbol\lesseqgtr 1351
\newsymbol\gtreqless 1352
\newsymbol\lesseqqgtr 1353
\newsymbol\gtreqqless 1354
\newsymbol\Rrightarrow 1356
\newsymbol\Lleftarrow 1357
\newsymbol\veebar 1259
\newsymbol\barwedge 125A
\newsymbol\doublebarwedge 125B
\undefine\angle
\newsymbol\angle 105C
\newsymbol\measuredangle 105D
\newsymbol\sphericalangle 105E
\newsymbol\varpropto 135F
\newsymbol\smallsmile 1360
\newsymbol\smallfrown 1361
\newsymbol\Subset 1362
\newsymbol\Supset 1363
\newsymbol\Cup 1264
  
\newsymbol\Cap 1265
  
\newsymbol\curlywedge 1266
\newsymbol\curlyvee 1267
\newsymbol\leftthreetimes 1268
\newsymbol\rightthreetimes 1269
\newsymbol\subseteqq 136A
\newsymbol\supseteqq 136B
\newsymbol\bumpeq 136C
\newsymbol\Bumpeq 136D
\newsymbol\lll 136E
  
\newsymbol\ggg 136F
  
\newsymbol\circledS 1073
\newsymbol\pitchfork 1374
\newsymbol\dotplus 1275
\newsymbol\backsim 1376
\newsymbol\backsimeq 1377
\newsymbol\complement 107B
\newsymbol\intercal 127C
\newsymbol\circledcirc 127D
\newsymbol\circledast 127E
\newsymbol\circleddash 127F
\newsymbol\lvertneqq 2300
\newsymbol\gvertneqq 2301
\newsymbol\nleq 2302
\newsymbol\ngeq 2303
\newsymbol\nless 2304
\newsymbol\ngtr 2305
\newsymbol\nprec 2306
\newsymbol\nsucc 2307
\newsymbol\lneqq 2308
\newsymbol\gneqq 2309
\newsymbol\nleqslant 230A
\newsymbol\ngeqslant 230B
\newsymbol\lneq 230C
\newsymbol\gneq 230D
\newsymbol\npreceq 230E
\newsymbol\nsucceq 230F
\newsymbol\precnsim 2310
\newsymbol\succnsim 2311
\newsymbol\lnsim 2312
\newsymbol\gnsim 2313
\newsymbol\nleqq 2314
\newsymbol\ngeqq 2315
\newsymbol\precneqq 2316
\newsymbol\succneqq 2317
\newsymbol\precnapprox 2318
\newsymbol\succnapprox 2319
\newsymbol\lnapprox 231A
\newsymbol\gnapprox 231B
\newsymbol\nsim 231C
\newsymbol\ncong 231D
\newsymbol\diagup 201E
\newsymbol\diagdown 201F
\newsymbol\varsubsetneq 2320
\newsymbol\varsupsetneq 2321
\newsymbol\nsubseteqq 2322
\newsymbol\nsupseteqq 2323
\newsymbol\subsetneqq 2324
\newsymbol\supsetneqq 2325
\newsymbol\varsubsetneqq 2326
\newsymbol\varsupsetneqq 2327
\newsymbol\subsetneq 2328
\newsymbol\supsetneq 2329
\newsymbol\nsubseteq 232A
\newsymbol\nsupseteq 232B
\newsymbol\nparallel 232C
\newsymbol\nmid 232D
\newsymbol\nshortmid 232E
\newsymbol\nshortparallel 232F
\newsymbol\nvdash 2330
\newsymbol\nVdash 2331
\newsymbol\nvDash 2332
\newsymbol\nVDash 2333
\newsymbol\ntrianglerighteq 2334
\newsymbol\ntrianglelefteq 2335
\newsymbol\ntriangleleft 2336
\newsymbol\ntriangleright 2337
\newsymbol\nleftarrow 2338
\newsymbol\nrightarrow 2339
\newsymbol\nLeftarrow 233A
\newsymbol\nRightarrow 233B
\newsymbol\nLeftrightarrow 233C
\newsymbol\nleftrightarrow 233D
\newsymbol\divideontimes 223E
\newsymbol\varnothing 203F
\newsymbol\nexists 2040
\newsymbol\Finv 2060
\newsymbol\Game 2061
\newsymbol\mho 2066
\newsymbol\eth 2067
\newsymbol\eqsim 2368
\newsymbol\beth 2069
\newsymbol\gimel 206A
\newsymbol\daleth 206B
\newsymbol\lessdot 236C
\newsymbol\gtrdot 236D
\newsymbol\ltimes 226E
\newsymbol\rtimes 226F
\newsymbol\shortmid 2370
\newsymbol\shortparallel 2371
\newsymbol\smallsetminus 2272
\newsymbol\thicksim 2373
\newsymbol\thickapprox 2374
\newsymbol\approxeq 2375
\newsymbol\succapprox 2376
\newsymbol\precapprox 2377
\newsymbol\curvearrowleft 2378
\newsymbol\curvearrowright 2379
\newsymbol\digamma 207A
\newsymbol\varkappa 207B
\newsymbol\Bbbk 207C
\newsymbol\hslash 207D
\undefine\hbar
\newsymbol\hbar 207E
\newsymbol\backepsilon 237F
\catcode`\@=\csname pre amssym.tex at\endcsname


\centerline{\twelvebf Sur le lemme fondamental pour les groupes
unitaires}
\vskip 5mm
\centerline{par G\'{e}rard Laumon\footnote{${}^{\ast}$}{\sevenrm CNRS
et Universit\'{e} Paris-Sud, UMR 8628, Math\'{e}matique, B\^{a}timent
425, F-91405 Orsay Cedex, France, Gerard.Laumon@math.u-psud.fr}}
\vskip 10mm

Dans [1], Goresky, Kottwitz et MacPherson formulent une conjecture de
puret\'{e} pour les fibres de Springer affines.  En admettant cette
conjecture, ils d\'{e}montrent un analogue g\'{e}om\'{e}trique du
lemme fondamental de Langlands et Shelstad pour un groupe r\'{e}ductif
$G$ sur un corps local non archim\'{e}dien $F$ d'\'{e}gales
caract\'{e}ristiques, en se limitant cependant aux \'{e}l\'{e}ments
r\'{e}guliers semi-simples de $G$ qui sont contenus dans une tore
maximal non ramifi\'{e}.

Dans cette note, nous d\'{e}duisons de cette m\^{e}me conjecture de
puret\'{e} pour les fibres de Springer affines de $\mathop{\rm GL}(n)$
un analogue g\'{e}om\'{e}trique du lemme fondamental pour les groupes
unitaires sur $F$, sans restriction sur l'\'{e}l\'{e}ment r\'{e}gulier
semi-simple.

Bien entendu, comme dans [1], il n'est pas difficile de d\'{e}duire de
cet \'{e}nonc\'{e} g\'{e}om\'{e}trique le lemme fondamental
arithm\'{e}tique pour les groupes unitaires sur $F$ par une
application directe de la formule des points fixes de
Grothendieck-Lefschetz (cf. [2] 1.3).

Notre approche est en partie inspir\'{e}e par celle de Goresky,
Kottwitz et MacPherson, en particulier en ce qui concerne l'usage de
la cohomologie \'{e}quivariante.  Elle en diff\`{e}re cependant sur un
point: nous faisons un usage essentiel du lien entre fibres de
Springer et jacobiennes compactifi\'{e}es d\'{e}velopp\'{e} dans [2].
\vskip 2mm

{\sevenrm Je remercie J.-B. Bost, M. Brion, L. Illusie, F. Loeser,
B.C. Ng\^{o}, M. Raynaud et J.-L. Waldspurger pour l'aide qu'ils m'ont
apport\'{e}e durant la pr\'{e}paration de ce travail.}

\section{1}{Lemme fondamental g\'{e}om\'{e}trique}

Soit $k$ un corps alg\'{e}briquement clos.  On consid\`{e}re une
famille finie $(E_{i})_{i\in I}$ d'extensions finies s\'{e}parables de
$F=k((\varpi_{F}))$, munies d'uniformisantes $\varpi_{E_{i}}\in {\cal
O}_{E_{i}}\subset E_{i}$, et un \'{e}l\'{e}ment
$$
\gamma_{I}=(\gamma_{i})_{i\in I}\in \bigoplus_{i\in I}\varpi_{E_{i}}
{\cal O}_{E_{i}}\subset {\cal O}_{E_{I}}=\bigoplus_{i\in I}{\cal
O}_{E_{i}}\subset E_{I}=\bigoplus_{i\in I}E_{i}
$$
tel que $E_{i}=F[\gamma_{i}]$ pour chaque $i\in I$ et que les
polyn\^{o}mes (unitaires) minimaux $P_{i}(x)\in {\cal O}_{F}[x]$ des
$\gamma_{i}$ soient deux \`{a} deux distincts.

Pour chaque partie $J$ de $I$, on appelle {\it fibre de Springer
affine en} $\gamma_{J}=(\gamma_{i})_{i\in J}$ et on note $X_{J}$ le
$k$-sch\'{e}ma des ${\cal O}_{F}$-r\'{e}seaux $M\subset
E_{J}$ tels que $\gamma_{J}M\subset M$.

Le groupe $\Lambda_{J}={\Bbb Z}^{J}$ agit sur $X_{J}$ par
$$
(\lambda_{i})_{i\in J}\cdot M=(\varpi_{E_{i}}^{-\lambda_{i}})_{i\in
J}M.
$$
Le groupe des composantes connexes de $X_{J}$ est le quotient
${\Bbb Z}^{J}\twoheadrightarrow {\Bbb Z},~\lambda\mapsto\sum_{i\in
J}\lambda_{i}$, et on a $X_{J}=\bigcup_{\lambda\in {\Bbb Z}}
X_{J}^{\lambda}={\Bbb Z}\times X_{J}^{0}$.  Le tore $T_{J}={\Bbb
G}_{{\rm m},k}^{J}$ agit sur $X_{J}$ via l'inclusion
$(k^{\times})^{J}\subset E_{J}^{\times}$.

Comme l'ont remarqu\'{e} Goresky, Kottwitz et MacPherson [1], la
conjecture suivante joue un r\^{o}le central dans l'\'{e}tude de la
cohomologie des fibres de Springer affines $X_{I}$.

\thm CONJECTURE 1.1 (Goresky, Kottwitz et MacPherson)
\enonce
Pour toute th\'{e}orie cohomologique de Weil, toute partie $J$ de $I$
et tout entier $n$, le $n$-\`{e}me groupe de cohomologie de $X_{J}$
est pur de poids $n$.
\endthm

Cette conjecture est \'{e}tablie dans les cas
{\og}{quasi-homog\`{e}nes}{\fg}.  Plus pr\'{e}cis\'{e}ment, Lusztig et
Smelt [3] l'ont d\'{e}montr\'{e}e dans le cas o\`{u} $I$ est
r\'{e}duit \`{a} un \'{e}l\'{e}ment $i$ et o\`{u} on peut choisir les
uniformisantes de telle sorte que $\varpi_{F}= \varpi_{E_{i}}^{n_{i}}$
et $\gamma_{i}= \varpi_{E_{i}}^{v_{i}}$ pour des entiers
$n_{i},v_{i}\geq 1$ premiers entre eux; en fait, sous ces
hypoth\`{e}ses, $X_{I}^{0}=X_{i}^{0}$ est projectif sur $k$ et
pav\'{e} par des espaces affines standard, de sorte que sa cohomologie
$\ell$-adique est m\^{e}me concentr\'{e}e en degr\'{e}s pairs.
Waldpurger, puis Goresky, Kottwitz et MacPherson, ont \'{e}tendu ce
r\'{e}sultat de Lusztig et Smelt au cas o\`{u}, pour des choix
convenables des uniformisantes, il existe des entiers $n,v\geq 1$
premiers entre eux et une famille $(\xi_{i})_{i\in I}$
d'\'{e}l\'{e}ments de $k$ deux \`{a} deux distincts tels que
$\varpi_{F}= \varpi_{E_{i}}^{n}$ et $\gamma_{i}=
\xi_{i}\varpi_{E_{i}}^{v}$ pour tout $i\in I$.
\vskip 3mm

Le quotient $Z_{J}$ de $X_{J}$ par l'action de $\Lambda_{J}$, ou ce
qui revient au m\^{e}me de $X_{J}^{0}$ par l'action du noyau
$\Lambda_{I}^{0}$ de l'\'{e}pimorphisme ${\Bbb Z}^{J}
\twoheadrightarrow {\Bbb Z}$ ci-dessus, est un $k$-sch\'{e}ma
projectif connexe de dimension
$$
\delta_{J}=\mathop{\rm dim}\nolimits_{k}({\cal O}_{E_{J}}/{\cal
O}_{F}[\gamma_{J}]).
$$

Fixons une partition non triviale $I=I_{1}\amalg I_{2}$, ou ce qui
revient au m\^{e}me le caract\`{e}re
$$
\kappa :\Lambda_{I}^{0}\rightarrow \{\pm 1\},~\lambda\mapsto
(-1)^{\sum_{i\in I_{1}}\lambda_{i}}= (-1)^{\sum_{i\in I_{2}}
\lambda_{i}}.
$$
Fixons aussi un nombre premier $\ell$ distinct de la
caract\'{e}ristique de $k$.

On a alors un syst\`{e}me local $\ell$-adique ${\cal L}$ de rang $1$
sur $Z_{I}$ d\'{e}fini par le rev\^{e}tement \'{e}tale galoisien
$X_{I}^{0}\rightarrow Z_{I}^{0}$ et le caract\`{e}re $\kappa$ de son
groupe de Galois $\Lambda_{I}^{0}$.  On a aussi une immersion
ferm\'{e}e
$$
X_{I_{1}}\times_{k}X_{I_{2}}\hookrightarrow X_{I}
$$
qui envoie le couple de ${\cal O}_{F}$-r\'{e}seaux $(M_{1}\subset
E_{I_{1}},M_{2}\subset E_{I_{2}})$ sur le ${\cal O}_{F}$-r\'{e}seau
$M=M_{1}\oplus M_{2}\subset E_{I}$.  Par passage au quotient par
$\Lambda_{I}=\Lambda_{I_{1}}\times\Lambda_{I_{2}}$, on en d\'{e}duit
une immersion ferm\'{e}e
$$
Z_{I_{1}}\times_{k}Z_{I_{2}}\hookrightarrow Z_{I}
$$
de codimension
$$
r=\delta_{I}-(\delta_{I_{1}}+\delta_{I_{2}})=\sum_{i\in
I_{1},j\in I_{2}}r_{ij}
$$
o\`{u} $r_{ij}=\mathop{\rm dim}\nolimits_{k}({\cal O}_{F}[x]/
(P_{i}(x),P_{j}(x))\geq 1$.

Si $Z_{I}$ et $Z_{I_{1}}\times_{k}Z_{I_{2}}$ \'{e}taient lisses sur
$k$, on aurait un morphisme de Gysin
$$
R\Gamma (Z_{I_{1}}\times_{k}Z_{I_{2}},{\Bbb Q}_{\ell})[-2r](-r)
\rightarrow R\Gamma (Z_{I},{\Bbb Q}_{\ell}),
$$
mais ce n'est pas du tout le cas.

La conjecture suivante, formul\'{e}e par Kottwitz, dit essentiellement
que, malgr\'{e} les singularit\'{e}s, on a quand m\^{e}me un morphisme
de Gysin pour cette immersion ferm\'{e}e, morphisme qui est en fait un
isomorphisme.  Le syst\`{e}me ${\cal L}$ semble y jouer le r\^{o}le d'un
faisceau d'orientation.

\thm CONJECTURE 1.2 (Lemme fondamental g\'{e}om\'{e}trique)
\enonce
Il existe un isomorphisme canonique
$$
H^{\bullet -2r}(Z_{I_{1}}\times_{k}Z_{I_{2}},{\Bbb Q}_{\ell})(-r)
\buildrel\sim\over\longrightarrow H^{\bullet}(Z_{I},{\cal L}).
$$
\endthm

On se propose d'esquisser dans cette note une d\'{e}monstration du
th\'{e}or\`{e}me suivant:

\thm TH\'{E}OR\`{E}ME 1.3
\enonce
La conjecture $1.1$ implique la conjecture $1.2$.
\endthm

\section{2}{Globalisation et d\'{e}formation {\rm (cf. [2])}}

Soit
$$
A=k[[\varpi_{F}]][x]/(P_{I}(x))
$$
o\`{u} pour toute partie $J$ de $I$, on a not\'{e}
$P_{J}(x)=\prod_{i\in J}P_{i}(x)$.  Alors $\mathop{\rm Spf}(A)$ est un
germe formel de courbe plane r\'{e}duite dont l'ensemble des branches
irr\'{e}ductibles est index\'{e} par $I$ (on rappelle que les
$P_{i}(x)\in (\varpi_{F},x)\subset k[[\varpi_{F}]][x]$ sont
irr\'{e}ductibles et deux \`{a} deux premiers entre eux).

On s'est donn\'{e} une partition non triviale $I=I_{1}\amalg I_{2}$.
On pose
$$
A_{\alpha}=k[[\varpi_{F}]][x]/(P_{I_{\alpha}}(x)),~\forall\alpha =1,2.
$$
L'entier $r$ du lemme fondamental est le nombre d'intersection des
germes formels de courbes $\mathop{\rm Spf}(A_{1})$ et $\mathop{\rm
Spf}(A_{2})$ trac\'{e}es sur le germe formel de surface $\mathop{\rm
Spf}(k[[\varpi_{F},x]])$, germes de courbes dont la r\'{e}union est
$\mathop{\rm Spf}(A)$.

On note $S=\mathop{\rm Spec}(k[[z]])$ {\og}{le}{\fg} trait strictement
hens\'{e}lien d'\'{e}gales caract\'{e}ristiques de corps r\'{e}siduel
$k$.  On note $s$ et $\eta$ les points ferm\'{e} et g\'{e}n\'{e}rique
de $S$, et on choisit un point g\'{e}om\'{e}trique $\overline{\eta}$
localis\'{e} en $\eta$ dont on note $K$ le corps r\'{e}siduel.

\thm PROPOSITION 2.1
\enonce
Il existe un morphisme de sch\'{e}mas $C\rightarrow S$ ayant les
propri\'{e}t\'{e}s suivantes:
\vskip 1mm

\itemitem{-} $C\rightarrow S$ est projectif et plat, \`{a} fibres
g\'{e}om\'{e}triques int\`{e}gres et de dimension $1$, de sorte que le
lieu singulier $D\subset C$ de $C\rightarrow S$ est fini sur $S$;
\vskip 1mm

\itemitem{-} la fibre sp\'{e}ciale r\'{e}duite de $D\rightarrow S$ est
constitu\'{e}e d'un seul point ferm\'{e} $c$ de $C_{s}$ et la fibre
g\'{e}n\'{e}rique g\'{e}om\'{e}trique r\'{e}duite de $D\rightarrow S$
est constitu\'{e}e de deux points ferm\'{e}s $c_{1},c_{2}$ rationnels
sur $\kappa (\eta )$ et de $r$ points ferm\'{e}s $(d_{ij,\rho})_{i\in
I_{1},j\in I_{2},\rho =1,\ldots ,r_{ij}}$;
\vskip 1mm

\itemitem{-} le compl\'{e}t\'{e} de l'anneau local de $C_{s}$ en
$c$ est isomorphe \`{a} $A$,
\vskip 1mm

\itemitem{-} le compl\'{e}t\'{e} de l'anneau local de
$C_{\overline{\eta}}$ en $c_{\alpha}$ est isomorphe \`{a}
$K\widehat{\otimes}_{k}A_{\alpha}$;
\vskip 1mm

\itemitem{-} la courbe $C_{\overline{\eta}}$ a une singularit\'{e}
quadratique ordinaire en chacun des points $d_{ij,\rho}$;
\vskip 1mm

\itemitem{-} $C$ admet une normalisation en famille par la droite
projective standard sur $S$,
$$
\nu :{\Bbb P}_{S}^{1}\rightarrow C
$$
telle que $\infty (S)\cap\nu^{-1}(D)=\emptyset$ o\`{u} $\infty :S
\rightarrow {\Bbb P}_{S}^{1}$ est la section infinie {\rm [} $\nu
:{\Bbb P}_{S}^{1}\rightarrow C$ est donc la normalisation de la
surface $S$, et les morphismes $\nu_{s}:{\Bbb P}_{s}^{1}\rightarrow
C_{s}$ et $\nu_{\overline{\eta}}: {\Bbb P}_{\overline{\eta}}^{1}
\rightarrow C_{\overline{\eta}}$ sont les normalisations des fibres
sp\'{e}ciale et g\'{e}n\'{e}rique g\'{e}om\'{e}trique de $C\rightarrow
S${\rm ]};
\vskip 1mm

\itemitem{-} le morphisme $\nu$ admet la factorisation
$$
{\Bbb P}_{S}^{1}\,\smash{\mathop{\hbox to 6mm{\rightarrowfill}}
\limits^{\scriptstyle \nu'}}\,C'\rightarrow C,
$$
en morphismes finis birationnels, o\`{u} $C_{s}'$ n'a comme seules
singularit\'{e}s que deux points $c_{1}'$ et $c_{2}'$ en lesquels les
compl\'{e}t\'{e}s des anneaux locaux de $C_{s}$ sont $A_{1}$ et
$A_{2}$, o\`{u} $\nu_{s}':{\Bbb P}_{k}^{1}\rightarrow C_{s}'$ est la
normalisation de $C_{s}'$, o\`{u}
$$
\nu'=S\times_{k}\nu_{s}:{\Bbb P}_{S}^{1}\rightarrow
C'=S\times_{k}C_{s}'
$$
et o\`{u} $C_{\overline{\eta}}'\rightarrow C_{\overline{\eta}}$ est la
normalisation des seules singularit\'{e}s quadratiques ordinaires
$d_{ij,\rho}$ de $C_{\overline{\eta}}$.
\endthm

\rem Preuve
\endrem
Nous nous contenterons ici de construire $C'\rightarrow C$ en laissant
au lecteur le soin de v\'{e}rifier les propri\'{e}t\'{e}s
annonc\'{e}es.

On construit tout d'abord comme dans la proposition 2.1.1 de [2] une
courbe $C_{s}'$ sur $k$ ayant les propri\'{e}t\'{e}s suivantes:
\vskip 1mm

\itemitem{-}  $C_{s}'$ est int\`{e}gre et projective sur $k$;
\vskip 1mm

\itemitem{-} $C_{s}'$ n'a comme seules singularit\'{e}s que deux
points ferm\'{e}s $c_{1,s}'$ et $c_{2,s}'$, points en lesquels les
compl\'{e}t\'{e}s des anneaux locaux de $C_{s}'$ sont respectivement
$A_{1}$ et $A_{2}$;
\vskip 1mm

\itemitem{-} $C_{s}'$ est normalis\'{e}e par la droite projective
standard ${\Bbb P}_{k}^{1}$, et le morphisme de normalisation
$\nu_{s}':{\Bbb P}_{k}^{1}\rightarrow C_{s}'$ envoie le point \`{a}
l'infini de ${\Bbb P}_{k}^{1}$ sur un point de $\infty_{s}\in C_{s}'$
qui est distinct de $c_{1,s}'$ et $c_{2,s}'$, et donc dans l'ouvert de
lissit\'{e} de $C_{s}'$.
\vskip 1mm

On consid\`{e}re la d\'{e}formation formelle constante
$$
{\Bbb P}_{\widehat{S}}^{1}\,\smash{\mathop{\hbox to 6mm
{\rightarrowfill}} \limits^{\scriptstyle \nu'}}
\,\widehat{C}'=\widehat{S}\times_{k}C_{s}'\rightarrow
\widehat{S}=\mathop{\rm Spf}(k[[z]])
$$
de $\nu_{s}':{\Bbb P}_{k}^{1}\rightarrow C_{s}'$.  Le
$\widehat{S}$-sch\'{e}ma formel $z$-adique $\widehat{C}'$ est
r\'{e}union de deux cartes affines, \`{a} savoir la restriction ${\cal
U}'=\mathop{\rm Spf}({\cal B}')$ de $\widehat{C}'$ \`{a} l'ouvert
$U'=C_{s}'-\{\infty_{s}\}$ de $C_{s}'=(\widehat{C}')_{{\rm red}}$ et
la restriction ${\cal V}'$ de $\widehat{C}'$ \`{a} l'ouvert
$V'=C_{s}'-\{c_{1,s}',c_{2,s}'\}\subset C_{s}'$.

On d\'{e}forme ensuite la normalisation partielle
$$\displaylines{
\qquad A=k[[\varpi_{F}]][x]/(P_{I}(x))
\hfill\cr\hfill
\subset k[[\varpi_{F}]][x]/(P_{I_{1}}(x))\times
k[[\varpi_{F}]][x]/(P_{I_{2}}(x))=A_{1}\times A_{2}=A'\qquad}
$$
en
$$\displaylines{
\qquad {\cal A}=k[[z,\varpi_{F}]][x]/(P_{I_{1}}(x))P_{I_{2}}(x-z))
\hfill\cr\hfill
\rightarrow k[[z,\varpi_{F}]][x]/(P_{I_{1}}(x))\times
k[[z,\varpi_{F}]][x](P_{I_{2}}(x))={\cal A}_{1}\times {\cal
A}_{2}={\cal A}'\qquad}
$$
o\`{u} cette derni\`{e}re fl\`{e}che envoie la classe de
$f(z,\varpi_{F},x)$ sur le couple form\'{e} de la classe de
$f(z,\varpi_{F},x)$ et de celle de $f(z,\varpi_{F},x+z)$.  Par
d\'{e}finition de $\widehat{C}'$, ${\cal A}_{\alpha}$ est le
compl\'{e}t\'{e} de l'anneau local en $c_{\alpha}'$ de $\widehat{C}'$,
ou ce qui revient au m\^{e}me de ${\cal U}'$.

On forme alors la $k[[z]]$-alg\`{e}bre $z$-adique produit fibr\'{e}
$$
{\cal B}={\cal A}\times_{{\cal A}'}{\cal B}'.
$$
La $k$-alg\`{e}bre r\'{e}duite $B=A\times_{A'}B'$ de ${\cal B}$ est
int\`{e}gre et de type finie sur $k$; la projection canonique
$B\hookrightarrow B'$ est un morphisme de normalisation partielle.

On v\'{e}rifie que ${\cal U}=\mathop{\rm Spf}({\cal B})$ est
formellement lisse sur $\widehat{S}$ en dehors d'un unique point
singulier $c_{s}\in U={\cal U}_{{\rm red}}=\mathop{\rm Spec}(B)$ et
que le compl\'{e}t\'{e} formel de l'anneau local de ${\cal U}$ en
$c_{s}$ n'est autre que ${\cal A}$.  Le morphisme ${\cal
U}'\rightarrow {\cal U}$ (d\'{e}fini par la projection canonique
${\cal B}\rightarrow {\cal B}'$) induit un isomorphisme de la
restriction de ${\cal U}'$ \`{a} l'ouvert
$U'-\{c_{1,s}',c_{2,s}'\}\subset U'$ sur la restriction de ${\cal U}$
\`{a} l'ouvert $U-\{c_{s}\}$ de $U$.

On construit alors une courbe formelle $z$-adique $\widehat{C}$ au
dessus de $\widehat{S}$ en recollant ${\cal U}$ et ${\cal V}'$ le long
de leur ouvert commun
$$
{\cal U}|_{U-\{c_{s}\}}\cong
{\cal U}'|_{U'-\{c_{1,s}',c_{2,s}'\}}={\cal V}'|_{V'-\{\infty_{s}\}},
$$
soit en d'autres termes comme la somme amalgam\'{e}e
$$
\widehat{C}=\mathop{\rm Spf}({\cal A})\amalg_{\mathop{\rm Spf}
({\cal A}')} \widehat{C}'.
$$

On d\'{e}finit enfin $C'\rightarrow C$ comme une alg\'{e}brisation du
morphisme canonique $\widehat{C}'\rightarrow \widehat{C}$.  La fibre
sp\'{e}ciale $C_{s}$ de $C\rightarrow S$ est la courbe obtenue en
recollant $U$ et $V'$ le long de leur ouvert commun $U-\{c_{s}\}\cong
U'-\{c_{1,s}',c_{2,s}'\}\cong V'-\{\infty_{s}\}$ et n'est autre qu'une
des courbes construites dans la proposition 2.1.1 de [2].

\hfill\hfill$\square$
\vskip 3mm

L'image inverse $\widetilde{D}$ de $D$ dans ${\Bbb P}_{S}^{1}$
contient une famille $(\widetilde{c}_{i})_{i\in I}$ de sections
constantes deux \`{a} deux distinctes telles que:
\vskip 1mm

\itemitem{-} $\widetilde{c}_{i}(s)$ corresponde \`{a} la branche
$\{P_{i}(x)=0\}$ de la singularit\'{e} de $C_{s}$ en $c$ pour
chaque $i\in I$,
\vskip 1mm

\itemitem{-} $\widetilde{c}_{i}(\overline{\eta})$ corresponde \`{a} la
branche $\{P_{i}(x)=0\}$ de la singularit\'{e} de
$C_{\overline{\eta}}$ en $c_{\alpha}$ pour chaque $i\in I_{\alpha}$
et $\alpha =1,2$.
\vskip 1mm

La fibre g\'{e}n\'{e}rique g\'{e}om\'{e}trique de $\widetilde{D}$
contient en outre une famille
$$
(\widetilde{d}_{ij,\rho}',\widetilde{d}_{ij\rho}'')_{i\in I_{1},j\in
I_{2},\rho =1,\ldots ,r_{ij}}
$$
de points ferm\'{e}s deux \`{a} deux distincts tels que:
\vskip 1mm

\itemitem{-} $\widetilde{d}_{ij,\rho}'$ et
$\widetilde{d}_{ij,\rho}''$ correspondent aux deux
branches de la singularit\'{e} quadratique ordinaire de
$C_{\overline{\eta}}$ en $d_{ij,\rho}$ pour tous $i\in I_{1}$, $j\in
I_{2}$ et $\rho =1,\ldots ,r_{ij}$,
\vskip 1mm

\itemitem{-} $\widetilde{d}_{ij,\rho}'$ et $\widetilde{d}_{ij,\rho}''$
se sp\'{e}cialisent en $\widetilde{c}_{i}(s)$ et
$\widetilde{c}_{j}(s)$ respectivement, quels que soient $i\in I_{1}$,
$j\in I_{2}$ et $\rho =1,\ldots ,r_{ij}$.
\vskip 3mm

Soient alors $P\rightarrow S$ et $g:\overline{P}\rightarrow S$ les
$S$-sch\'{e}mas de Picard et de Picard compactifi\'{e} relatifs de
$C\rightarrow S$.  On sait que $P={\Bbb Z}\times P^{0}$ o\`{u} $P^{0}$
est un $S$-sch\'{e}ma affine lisse en groupes commutatifs, que
$\overline{P}={\Bbb Z}\times\overline{P}{}^{\,0}$ o\`{u}
$\overline{P}{}^{\,0}$ est un $S$-sch\'{e}ma projectif et plat, \`{a}
fibres g\'{e}om\'{e}triquement int\`{e}gres et localement
d'intersection compl\`{e}te, qui contient $P^{0}$ comme un ouvert
dense, et que l'action par translation de $P$ sur lui-m\^{e}me se
prolonge en une action de $P$ sur $\overline{P}$ compatible en un sens
\'{e}vident avec l'action par translation de ${\Bbb Z}$ sur
lui-m\^{e}me.

Le tore maximal $T_{s}$ de la jacobienne $P_{s}^{0}$ de $C_{s}$ est
canoniquement isomorphe \`{a} ${\Bbb G}_{{\rm m},s}^{I}/{\Bbb
G}_{{\rm m},s}$ et le tore maximal $\widetilde{T}_{\overline{\eta}}$
de la jacobienne $P_{\overline{\eta}}^{0}$ de $C_{\overline{\eta}}$ est
canoniquement isomorphe \`{a}
$$
({\Bbb G}_{{\rm m},K}^{I_{1}}/{\Bbb G}_{{\rm m},K}) \times_{K}
({\Bbb G}_{{\rm m},K}^{I_{1}}/{\Bbb G}_{{\rm m},K})\times_{K}
\prod_{{\scriptstyle i\in I_{1}\atop\scriptstyle j\in I_{2}}}
({\Bbb G}_{{\rm m},K}^{2}/{\Bbb G}_{{\rm m},K})^{r_{ij}}.
$$
De plus, le tore maximal $T_{s}\subset P_{s}^{0}$ se rel\`{e}ve en un sous-tore
$$
T={\Bbb G}_{{\rm m},S}^{I}/{\Bbb G}_{{\rm m},S}\subset P^{0}
$$
dont la fibre g\'{e}n\'{e}rale
$$
T_{\overline{\eta}}={\Bbb G}_{{\rm m},K}^{I}/{\Bbb G}_{{\rm m},
\overline{\eta}}\subset ({\Bbb G}_{{\rm m},K}^{I_{1}}/{\Bbb
G}_{{\rm m},K}) \times_{K}({\Bbb G}_{{\rm m},K}^{I_{1}}/{\Bbb G}_{{\rm
m},K})\times_{K}\prod_{{\scriptstyle i\in I_{1}\atop\scriptstyle j\in
I_{2}}}({\Bbb G}_{{\rm m},K}^{2}/{\Bbb G}_{{\rm m},K})^{r_{ij}}
=\widetilde{T}_{\overline{\eta}}
$$
est d\'{e}finie par l'inclusion
$$
(u_{i})_{i\in I}\mapsto ((v_{1,i})_{i\in I_{1}},(v_{2,i})_{i\in
I_{2}},(w_{ij,\rho}',w_{ij,\rho}'')_{i\in I_{1},j\in I_{2},\rho
=1,\ldots ,r_{ij}})
$$
o\`{u} $v_{\alpha ,i}=u_{i}$ pour tout $i\in I_{\alpha}$ et $\alpha
=1,2$, et o\`{u} $w_{ij,\rho}'=u_{i}$ et $w_{ij,\rho}''=u_{j}$ pour
tous $i\in I_{1}$, $j\in I_{2}$ et $\rho =1,\ldots ,r_{ij}$. (On
rappelle que les $r_{ij}$ sont $\geq 1$.)
\vskip 3mm

Par auto-dualit\'{e} (partielle) des jacobiennes compactifi\'{e}es, on a:
\vskip 1mm

\itemitem{-} un rev\^{e}tement \'{e}tale galoisien
$T_{s}$-\'{e}quivariant $\pi_{s}:\overline{P}{}_{s}^{\,\natural}
\rightarrow \overline{P}_{s}$ de groupe de Galois le groupe
$X^{\ast}(T_{s})$ des caract\`{e}res de $T_{s}$, c'est-\`{a}-dire le
noyau du morphisme somme ${\Bbb Z}^{I}\rightarrow {\Bbb Z}$,
\vskip 1mm

\itemitem{-} un rev\^{e}tement \'{e}tale galoisien
$\widetilde{T}_{\overline{\eta}}$-\'{e}quivariant
$\widetilde{\pi}_{\overline{\eta}}: \widetilde{\overline{P}}
{}_{\overline{\eta}}^{\,\natural}\rightarrow
\overline{P}_{\overline{\eta}}$ de groupe de Galois le groupe
$X^{\ast}(\widetilde{T}_{\overline{\eta}})$ des caract\`{e}res de
$\widetilde{T}_{\overline{\eta}}$, c'est-\`{a}-dire
$$
X^{\ast}(\widetilde{T}_{\overline{\eta}})=\mathop{\rm Ker}
({\Bbb Z}^{I_{1}}\rightarrow {\Bbb  Z})\times \mathop{\rm Ker}
({\Bbb Z}^{I_{2}}\rightarrow {\Bbb Z})\times \prod_{{\scriptstyle
i\in I_{1}\atop\scriptstyle j\in I_{2}}}(\mathop{\rm Ker}
({\Bbb Z}^{2}\rightarrow {\Bbb Z}))^{r_{ij}},
$$
\vskip 1mm

\itemitem{-} un rev\^{e}tement \'{e}tale galoisien $T$-\'{e}quivariant
$\pi :\overline{P}{}^{\,\natural}\rightarrow \overline{P}$ de groupe de
Galois $X^{\ast}(T)=X^{\ast}(T_{s})$, dont la fibre sp\'{e}ciale est
$\pi_{s}:\overline{P}{}_{s}^{\,\natural}\rightarrow \overline{P}_{s}$ et
dont la fibre g\'{e}n\'{e}rale $\pi_{\overline{\eta}}:
\overline{P}{}_{\overline{\eta}}^{\,\natural} \rightarrow
\overline{P}_{\overline{\eta}}$ est le quotient du rev\^{e}tement
$\widetilde{\pi}_{\overline{\eta}}:
\widetilde{\overline{P}}{}_{\overline{\eta}}^{\,\natural} \rightarrow
\overline{P}_{\overline{\eta}}$ correspondant au quotient
$$\displaylines{
\qquad X^{\ast}(\widetilde{T}_{\overline{\eta}})=\mathop{\rm Ker}
({\Bbb Z}^{I_{1}}\rightarrow {\Bbb Z})\times \mathop{\rm Ker}
({\Bbb Z}^{I_{2}}\rightarrow {\Bbb Z})\times\prod_{{\scriptstyle
i\in I_{1}\atop\scriptstyle j\in I_{2}}}(\mathop{\rm Ker}
({\Bbb Z}^{2}\rightarrow {\Bbb Z}))^{r_{ij}}
\hfill\cr\hfill
\twoheadrightarrow \mathop{\rm Ker} ({\Bbb Z}^{I}\rightarrow {\Bbb
Z})=X^{\ast}(T_{\overline{\eta}})=X^{\ast}(T)\qquad}
$$
du groupe de Galois, quotient qui s'\'{e}crit concr\`{e}tement
$$
\bigl((\lambda_{1,i})_{i\in I_{1}},(\lambda_{2,i})_{i\in
I_{2}},(\mu_{ij,\rho})_{i\in I_{1},j\in I_{2},\rho =1,\ldots
,r_{ij}}\bigr)\mapsto (\lambda_{i})_{i\in I},
$$
o\`{u}
$$
\lambda_{i}=\lambda_{1,i}+\sum_{j\in I_{2}}\sum_{\rho =1}^{r_{ij}}
\mu_{ij,\rho}',~\forall i\in I_{1},\hbox{ et }\lambda_{j}=
\lambda_{2,j}+\sum_{i\in I_{1}}\sum_{\rho =1}^{r_{ij}}
\mu_{ij,\rho}'',~\forall j\in I_{2}.
$$
\vskip 3mm

On note $V$ le $k$-sch\'{e}ma localement de type fini
obtenu en prenant ${\Bbb Z}$ copies de la droite projective standard
${\Bbb P}_{k}^{1}$ sur $k$ et en identifiant le point \`{a} l'infini
de la $n$-\`{e}me copie avec l'origine de la $(n+1)$-\`{e}me.  Le
groupe ${\Bbb Z}$ agit proprement et librement par translation sur
$V$ et le quotient de $V$ par cette action est
la cubique plane \`{a} point double ordinaire $W$ obtenue en
identifiant l'origine et le point \`{a} l'infini dans ${\Bbb
P}_{k}^{1}$.  On v\'{e}rifie que ${\Bbb Z}\times V$ est la
fibre de Springer affine form\'{e}e des ${\cal O}_{F}$-r\'{e}seaux
$N\subset F\oplus F$ tels que $(\varpi_{F},-\varpi_{F})N\subset N$.

Le rev\^{e}tement \'{e}tale galoisien $V\rightarrow
W=V/{\Bbb Z}$ de groupe de Galois ${\Bbb Z}$ est ${\Bbb
G}_{{\rm m},k}$-\'{e}quivariant pour les actions induites par l'action
standard de ${\Bbb G}_{{\rm m},k}$ sur ${\Bbb P}_{k}^{1}$. On a
$$
V^{{\Bbb G}_{{\rm m},k}}={\Bbb Z}\times\mathop{\rm
Spec}(k)
$$
et
$$
W^{{\Bbb G}_{{\rm m},k}}=\mathop{\rm Spec}(k)
$$
(les points fixes sont les points doubles de $V$ et $W$).

\thm LEMME 2.2
\enonce
Le rev\^{e}tement \'{e}tale galoisien $\overline{P}{}_{s}^{\,\natural
,0}\rightarrow \overline{P}{}_{s}^{\,0}$ de groupe de Galois
$X^{\ast}(T)$ est hom\'{e}omorphe de mani\`{e}re \'{e}quivariante sous
$T_{s}\cong T_{I}/{\Bbb G}_{{\rm m},k}$ au rev\^{e}tement \'{e}tale
galoisien $X_{I}^{0}\rightarrow Z_{I}=X_{I}^{0}/\Lambda_{I}^{0}$ de
groupe de Galois $\Lambda_{I}^{0}\cong X^{\ast}(T)$.

Le rev\^{e}tement $\widetilde{\pi}_{\overline{\eta}}:
\widetilde{\overline{P}}{}_{\overline{\eta}}^{\,\natural ,0}\rightarrow
\overline{P}{}_{\overline{\eta}}^{\,0}$ de groupe de Galois
$X^{\ast}(\widetilde{T}_{\overline{\eta}})$ est hom\'{e}omorphe de
mani\`{e}re \'{e}quivariante sous
$$
\widetilde{T}_{\overline{\eta}}\cong K\otimes_{k}((T_{I_{1}}/
{\Bbb G}_{{\rm m},k})\times_{k}(T_{I_{2}}/{\Bbb G}_{{\rm m},k})
\times_{k}{\Bbb G}_{{\rm m},k}^{r})
$$
au rev\^{e}tement
$$
K\otimes_{k}\bigl(X_{I_{1}}^{0}\times_{k}X_{I_{2}}^{0}\times_{k}
V^{r}\rightarrow Z_{I_{1}}\times_{k} Z_{I_{2}}
\times_{k}W^{r}\bigr)
$$
de groupe de Galois $\Lambda_{I_{1}}^{0}\times\Lambda_{I_{2}}^{0}
\times {\Bbb Z}^{r}\cong X^{\ast}(\widetilde{T}_{\overline{\eta}})$.
\endthm

\rem Preuve
\endrem
On rappelle qu'il revient au m\^{e}me de se donner un ${\cal
O}_{C_{s}}$-Module coh\'{e}rent sans torsion de rang g\'{e}n\'{e}rique
$1$, muni d'une trivialisation sur $C_{s}-\{c\}$, ou de se donner un
${\cal O}_{F}$-r\'{e}seau $M\subset E_{I}$ tel que $\gamma_{I}M\subset
M$.

De m\^{e}me, il est \'{e}quivalent de se donner un ${\cal
O}_{C_{\overline{\eta}}}$-Module coh\'{e}rent sans torsion de rang
g\'{e}n\'{e}rique $1$, muni d'une trivialisation sur
$C_{\overline{\eta}}-\{c_{1},c_{2},(d_{ij,\rho})_{i\in I_{1},j\in
I_{2},\rho =1,\ldots ,r_{ij}}\}$, ou de se donner un uplet de ${\cal
O}_{F}$-r\'{e}seaux
$$
(M_{1}\subset E_{I_{1}},M_{2}\subset E_{I_{2}},(N_{ij,\rho}\subset
F\oplus F)_{i\in I_{1},j\in I_{2},\rho =1,\ldots ,r_{ij}})
$$
tels que $\gamma_{I_{1}}M_{1}\subset M_{1}$,
$\gamma_{I_{2}}M_{2}\subset M_{2}$ et
$(\varpi_{F},-\varpi_{F})N_{ij,\rho}\subset N_{ij,\rho}$.

\hfill\hfill$\square$
\vskip 3mm

Restreignons nous \`{a} partir de maintenant \`{a} l'action du sous-tore
$$
{\Bbb G}_{{\rm m},S}\buildrel\sim\over\longrightarrow {\Bbb G}_{{\rm
m},S}^{2}/{\Bbb G}_{{\rm m},S}\hookrightarrow {\Bbb G}_{{\rm m},S}^{I}
/{\Bbb G}_{{\rm m},S}=T
$$
(et donc du sous-tore ${\Bbb G}_{{\rm m},K}\subset
T_{\overline{\eta}}\subset \widetilde{T}_{\overline{\eta}}$), o\`{u}
l'isomorphisme de gauche envoie $u$ sur la classe de $(u,1)$ et o\`{u}
le monomorphisme de droite envoie la classe $(u_{1},u_{2})$ sur celle
de $(t_{i})_{i\in I}$ avec $t_{i}=u_{\alpha}$ quel que soient $i\in
I_{\alpha}$ et $\alpha =1,2$.

On a donc:
\vskip 1mm

\itemitem{-} une action de ${\Bbb G}_{{\rm m},S}$ sur $\overline{P}$,
\vskip 1mm

\itemitem{-} un rev\^{e}tement ${\Bbb G}_{{\rm m},S}$-\'{e}quivariant
$\pi :\overline{P}{}^{\,\natural}\rightarrow \overline{P}$ de groupe de
Galois $X^{\ast}(T)$,
\vskip 1mm

\itemitem{-} un rev\^{e}tement ${\Bbb G}_{{\rm m},K}$-\'{e}quivariant
$\widetilde{\pi}_{\overline{\eta}}:\widetilde{\overline{P}}
{}_{\overline{\eta}}^{\,\natural}\rightarrow
\overline{P}_{\overline{\eta}}$ de groupe de Galois
$X^{\ast}(\widetilde{T}_{\overline{\eta}})$,
\vskip 1mm

\itemitem{-} une identification de $\pi_{\overline{\eta}}$ comme le
quotient de $\widetilde{\pi}_{\overline{\eta}}$ correspondant \`{a}
un \'{e}pimorphisme $X^{\ast}(\widetilde{T}_{\overline{\eta}})
\twoheadrightarrow X^{\ast}(T_{\overline{\eta}})$.
\vskip 3mm

\thm LEMME 2.3
\enonce
Le lieu des points fixes $(\overline{P}{}^{\,0})^{{\Bbb G}_{{\rm
m},S}}\subset\overline{P}{}^{\,0}$ pour l'action du sous-tore ${\Bbb
G}_{{\rm m},S}$ ci-dessus est $S$-hom\'{e}omorphe \`{a}
$$
S\times_{k}(Z_{I_{1}}\times_{k}Z_{I_{2}})
$$
et les immersion ferm\'{e}es $(\overline{P}{}_{s}^{\,0})^{{\Bbb
G}_{{\rm m},s}}\hookrightarrow \overline{P}{}_{s}^{\,0}$ et
$(\overline{P} {}_{\overline{\eta}}^{\,0})^{{\Bbb G}_{{\rm
m},\overline{\eta}}} \hookrightarrow
\overline{P}{}_{\overline{\eta}}^{\,0}$ correspondantes sont
hom\'{e}omorphes aux immersions ferm\'{e}es
$$
Z_{I_{1}}\times_{k}Z_{I_{2}}\hookrightarrow Z_{I}\hbox{ et }
K\otimes_{k}(Z_{I_{1}}\times_{k}Z_{I_{2}}\times_{k}(\mathop{\rm
Spec}(k)\hookrightarrow W^{r}))
$$
respectivement, o\`{u} $\mathop{\rm Spec}(k)\hookrightarrow W^{r}$
est le $r$-uplet dont toutes les composantes sont le point double
de $W$.

Le lieu des points fixes $(\overline{P}{}^{\,\natural ,0})^{{\Bbb
G}_{{\rm m},S}}\subset \overline{P}{}^{\,\natural ,0}$ pour l'action
du sous-tore ${\Bbb G}_{{\rm m},S}$ ci-dessus est $S$-hom\'{e}omorphe
\`{a}
$$
S\times_{k}\bigcup_{\lambda\in {\Bbb Z}}(X_{I_{1}}^{\lambda}
\times_{k}X_{I_{2}}^{-\lambda})= S\times_{k}({\Bbb Z}\times
(X_{I_{1}}^{0} \times_{k}X_{I_{2}}^{0})).
$$
et les immersion ferm\'{e}es $(\overline{P}{}_{s}^{\,\natural
,0})^{{\Bbb G}_{{\rm m},s}}\hookrightarrow \overline{P}
{}_{s}^{\,\natural ,0}$ et $(\overline{P}
{}_{\overline{\eta}}^{\,\natural ,0})^{{\Bbb G}_{{\rm m},
\overline{\eta}}} \hookrightarrow
\overline{P}{}_{\overline{\eta}}^{\,\natural ,0}$ correspondantes sont
hom\'{e}omorphes aux immersions ferm\'{e}es
$$
{\Bbb Z}\times (X_{I_{1}}^{0}\times_{k}X_{I_{2}}^{0})\hookrightarrow
X_{I}^{0}
$$
et
$$
K\otimes_{k}(X_{I_{1}}^{0} \times_{k}X_{I_{2}}^{0}\times_{k}({\Bbb Z}
\times \mathop{\rm Spec}(k)\hookrightarrow (V^{r}/\mathop{\rm
Ker}({\Bbb Z}^{r}\rightarrow {\Bbb Z}))))
$$
respectivement, o\`{u} $V^{r}/\mathop{\rm Ker}({\Bbb Z}^{r}\rightarrow
{\Bbb Z})$ est le rev\^{e}tement \'{e}tale galoisien de $W^{r}$ de
groupe de Galois ${\Bbb Z}$ qui est le quotient du rev\^{e}tement
\'{e}tale galoisien $V^{r}\rightarrow W^{r}$ correspondant au quotient
${\Bbb Z}^{r} \twoheadrightarrow {\Bbb Z},~(\lambda_{\rho})_{\rho
=1,\ldots ,r} \mapsto \sum_{\rho =1}^{r}\lambda_{\rho}$, de son groupe
de Galois, et o\`{u} ${\Bbb Z}\times\mathop{\rm Spec}(k)
\hookrightarrow (V^{r}/\mathop{\rm Ker}({\Bbb Z}^{r} \rightarrow {\Bbb
Z}))$ est le lieu des points fixes sous l'action diagonale de ${\Bbb
G}_{{\rm m},k}$.

L'immersion ferm\'{e}e $(\widetilde{\overline{P}}
{}_{\overline{\eta}}^{\,\natural ,0})^{{\Bbb G}_{{\rm m},
\overline{\eta}}}\subset \widetilde{\overline{P}}
{}_{\overline{\eta}}^{\,\natural ,0}$ du lieu des points fixes pour
l'action du sous-tore ${\Bbb G}_{{\rm m},\overline{\eta}}$ ci-dessus
est hom\'{e}omorphe \`{a}
$$
K\otimes_{k}(X_{I_{1}}^{0}\times_{k}X_{I_{2}}^{0}
\times_{k}({\Bbb Z}^{r}\times \mathop{\rm Spec}(k)\hookrightarrow
V^{r}))
$$
o\`{u} ${\Bbb Z}^{r}\times\mathop{\rm Spec}(k) \hookrightarrow V^{r}$
est le lieu des points fixes sous l'action diagonale de ${\Bbb
G}_{{\rm m},k}$.

\hfill\hfill$\square$
\endthm

\rem Preuve
\endrem
Le ferm\'{e} des points fixes $\overline{P}{}^{\,{\Bbb G}_{{\rm m},S}}
\hookrightarrow \overline{P}$ admet la description suivante.  On
consid\`{e}re la normalisation partielle interm\'{e}diaire
$$
{\Bbb P}_{S}^{1}\buildrel\nu'\over\rightarrow C'\rightarrow C
$$
de la proposition $2.1$.  Alors, l'op\'{e}ration d'image directe pour
les ${\cal O}_{C'}$-Modules sans torsion de rang g\'{e}n\'{e}rique $1$
par le morphisme fini birationnel $C'\rightarrow C$ d\'{e}finit un
$S$-morphisme $\iota :\overline{P}{}'\rightarrow \overline{P}{}$ du
$S$-sch\'{e}ma de Picard compactifi\'{e} $\overline{P}{}'$ de $C'$ dans
$\overline{P}$, et l'image de $\iota$ est pr\'{e}cis\'{e}ment
$\overline{P}^{{\Bbb G}_{{\rm m},S}}$.  De plus, $\iota$ est un
hom\'{e}omorphisme sur son image.

Maintenant $\overline{P}{}'$ est isomorphe \`{a} $S\times_{k}
\overline{P}{}_{s}'$, $\iota_{s}:(\overline{P}{}_{s}')^{0}\rightarrow
\overline{P}{}_{s}^{\,0}$ est hom\'{e}omorphe \`{a}
$Z_{I_{1}}\times_{k} Z_{I_{2}}\hookrightarrow Z_{I}$ et
$\iota_{\overline{\eta}}:(\overline{P}{}_{\overline{\eta}}')^{0}\rightarrow
\overline{P}{}_{\overline{\eta}}^{\,0}$ est hom\'{e}omorphe \`{a}
$K\otimes_{k}(Z_{I_{1}}\times_{k}Z_{I_{2}}\times_{k}(\mathop{\rm
Spec}(k)\hookrightarrow W^{r}))$.

\hfill\hfill$\square$
\vskip 3mm

\section{3}{Sp\'{e}cialisation en homologie \'{e}quivariante}

Consid\'{e}rons un $S$-sch\'{e}ma $f:X\rightarrow S$ muni d'une action
du tore $T={\Bbb G}_{{\rm m},k}^{m}$ et d'une action libre d'un ${\Bbb
Z}$-module libre de rang fini $\Lambda$, qui commute \`{a} celle de
$T$, $X$ \'{e}tant r\'{e}union filtrante croissante d'une famille de
ferm\'{e}s $(X_{n})_{n}$ qui sont propres sur $S$ et stables sous $T$.
En pratique, $X$ sera notre $\overline{P}{}^{\,\natural}$ muni de
l'action de $X^{\ast}(T)$ et, si n\'{e}cessaire, on se limitera \`{a}
ce cas.

Pour travailler avec des modules du type ${\Bbb
Q}_{\ell}[\Lambda ]$ plut\^{o}t qu'avec des modules du type ${\Bbb
Q}_{\ell}^{\Lambda}$, on {\og}{dualise}{\fg} la cohomologie $\ell$-adique
$T$-\'{e}quivariante et on utilise donc l'homologie $\ell$-adique
$T$-\'{e}quivariante.

On rappelle que l'homologie \'{e}quivariante
$$
H_{\bullet}^{T}(\mathop{\rm Spec}(k),{\Bbb Q}_{\ell})
$$
du point est un module gradu\'{e} sur la cohomologie \'{e}quivariante
$H_{T}^{\bullet}(\mathop{\rm Spec}(k),{\Bbb Q}_{\ell})\cong {\Bbb
Q}_{\ell}[\delta_{1},\ldots ,\delta_{m}]$, module qui est isomorphe
\`{a} ${\Bbb Q}_{\ell}[t_{1},\ldots ,t_{m}]$ avec $\mathop{\rm
deg}(t_{\mu})=2$ sur lequel les $\delta_{\mu}$ agissent comme des
d\'{e}rivations normalis\'{e}es par $\delta_{\mu}
(t_{\mu'})=\delta_{\mu ,\mu'}$, et donc sur lequel ${\Bbb
Q}_{\ell}[\delta_{1},\ldots ,\delta_{m}]$ agit comme un anneau
d'op\'{e}rateurs diff\'{e}rentiels \`{a} coefficients constants.

On a des ${\Bbb Q}_{\ell}[\Lambda ][\delta_{1},\ldots
,\delta_{m}]$-modules
$$
H_{\bullet}^{T}(X_{s},{\Bbb Q}_{\ell})=\lim_{{\scriptstyle
\longrightarrow\atop\scriptstyle n}}\mathop{\rm Hom}\nolimits_{{\Bbb
Q}_{\ell}}(H_{T}^{\bullet}(X_{n,s},{\Bbb Q}_{\ell}),{\Bbb Q}_{\ell})
$$
et
$$
H_{\bullet}^{T}(X_{\overline{\eta}},{\Bbb Q}_{\ell})=
\lim_{{\scriptstyle \longrightarrow\atop\scriptstyle n}}\mathop{\rm
Hom}\nolimits_{{\Bbb Q}_{\ell}}(H_{T}^{\bullet}
(X_{n,\overline{\eta}},{\Bbb Q}_{\ell}),{\Bbb Q}_{\ell}).
$$

\thm PROPOSITION 3.1
\enonce
On a un carr\'{e} commutatif canonique
$$\diagram{
H_{\bullet}^{T} (X_{\overline{\eta}}^{T},{\Bbb Q}_{\ell}) &\kern -1mm
\smash{\mathop{\hbox to 8mm{\rightarrowfill}}\limits^{\scriptstyle}}
\kern -1mm&H_{\bullet}^{T}(X_{s}^{T}, {\Bbb Q}_{\ell})\cr
\llap{$\scriptstyle $}\left\downarrow
\vbox to 4mm{}\right.\rlap{}&&\llap{}\left\downarrow
\vbox to 4mm{}\right.\rlap{$\scriptstyle $}\cr
H_{\bullet}^{T}(X_{\overline{\eta}}, {\Bbb Q}_{\ell})&\kern -1mm
\smash{\mathop{\hbox to 8mm{\rightarrowfill}}\limits^{\scriptstyle }}
\kern -1mm& H_{\bullet}^{T}(X_{s},{\Bbb Q}_{\ell})\cr}
$$
de ${\Bbb Q}_{\ell}[\Lambda][\delta_{1},\ldots ,\delta_{m}]$-modules
gradu\'{e}s, o\`{u} les fl\`{e}ches horizontales sont les fl\`{e}ches
de sp\'{e}cialisation pour le $S$-sch\'{e}ma des points fixes $X^{T}$ et
$X$, et o\`{u} les fl\`{e}ches verticales sont les fl\`{e}ches de 
restriction au
ferm\'{e} $X^{T}\subset X$.
\endthm

Pour d\'{e}montrer la proposition, nous aurons besoin du lemme suivant:

\thm LEMME 3.2
\enonce
Notons $[S/T]$ le $S$-champ alg\'{e}brique quotient de $S$ par
l'action de $T$ et $\varepsilon :[S/T]\rightarrow S$ son morphisme
structural.

Pour tout objet $K$ de $D_{{\rm c}}^{{\rm b}}([S/T],{\Bbb Q}_{\ell})$, la
formation de $R\varepsilon_{\ast}K$ commute a tout changement de base
$S'\rightarrow S$.

En particulier, on a une fl\`{e}che de sp\'{e}cialisation en cohomologie
$\ell$-adique $T$-\'{e}qui\-va\-riante
$$
R\Gamma_{T}(s,K_{s})\rightarrow R\Gamma_{T}(\overline{\eta},
K_{\overline{\eta}}).
$$
\endthm

\rem Preuve
\endrem
Il suffit de d\'{e}montrer le lemme pour tout objet $K$ de
$D_{{\rm c}}^{[a,b]}([S/T],{\Bbb Q}_{\ell})$ quel que soient les
entiers $a\leq b$.

Pour simplifier, nous nous limiterons au cas o\`{u} $T={\Bbb G}_{{\rm
m},k}$.  Soient $p:{\Bbb P}_{S}^{N}\rightarrow [S/T]$ le morphisme
repr\'{e}sentable induit par passage au quotient du morphisme
structural ${\Bbb A}_{S}^{N+1}-0(S)\rightarrow S$, o\`{u}
$0:S\rightarrow {\Bbb A}_{S}^{N+1}$ est la section nulle.  Le
morphisme $p$ est universellement $2N$-acyclique.  Par suite, si $K\in
D_{{\rm c}}^{[a,b]}([S/T],{\Bbb Q}_{\ell})$ et si $2N\geq b-a$, la
fl\`{e}che d'adjonction $K\rightarrow Rp_{\ast}p^{\ast}K$ induit un
isomorphisme
$$
K\buildrel\sim\over\longrightarrow \tau_{\leq b}Rp_{\ast}p^{\ast}K,
$$
et donc un isomorphisme
$$
\tau_{\leq b}R\varepsilon_{\ast}K\buildrel\sim\over\longrightarrow
\tau_{\leq b}R\varepsilon_{\ast}\tau_{\leq b}Rp_{\ast}p^{\ast}K=
\tau_{\leq b}R\varepsilon_{\ast}Rp_{\ast}p^{\ast}K.
$$

Or la projection canonique $\varepsilon\circ p:{\Bbb P}_{S}^{N}
\rightarrow S$ est propre et le morphisme $p$ est lisse, de sorte que
la formation de $R\varepsilon_{\ast}Rp_{\ast}p^{\ast}K$ commute \`{a}
tout changement de base.
\hfill\hfill$\square$
\vskip 3mm

\rem Preuve de la proposition
\endrem
Comme chaque $f_{n}=f|X_{n}:X_{n}\rightarrow S$ est un $S$-sch\'{e}ma
propre, la formation de $Rf_{n,\ast}{\Bbb Q}_{\ell ,X}$ commute \`{a}
tout changement de base $S'\rightarrow S$ et on peut appliquer le
lemme 3.2 \`{a} $K=R\overline{f}_{n,\ast}{\Bbb Q}_{\ell ,[X_{n}/T]}$
o\`{u} on a not\'{e}
$$
\overline{f}_{n}:[X_{n}/T]\rightarrow [S/T]
$$
le morphisme repr\'{e}sentable et propre induit par $f_{n}$.  En
particulier, on a une fl\`{e}che de sp\'{e}cialisation en cohomologie
$\ell$-adique $T$-\'{e}quivariante
$$
R\Gamma_{T}(X_{n,s},{\Bbb Q}_{\ell})\rightarrow
R\Gamma_{T}(X_{n,\overline{\eta}},{\Bbb Q}_{\ell}).
$$
et un carr\'{e} commutatif de ${\Bbb Q}_{\ell}
[\Lambda][\delta_{1},\ldots ,\delta_{m}]$-modules gradu\'{e}s
$$\diagram{
H_{T}^{\bullet}(X_{n,s},{\Bbb Q}_{\ell})&\kern -1mm\smash{\mathop{\hbox
to 8mm{\rightarrowfill}} \limits^{\scriptstyle }}\kern -1mm&
H_{T}^{\bullet}(X_{n,\overline{\eta}},{\Bbb Q}_{\ell})\cr
\llap{$\scriptstyle $}\left\downarrow
\vbox to 4mm{}\right.\rlap{}&&\llap{}\left\downarrow
\vbox to 4mm{}\right.\rlap{$\scriptstyle $}\cr
H_{T}^{\bullet}(X_{n,s}^{T},{\Bbb Q}_{\ell})&\kern -1mm
\smash{\mathop{\hbox to 8mm{\rightarrowfill}}
\limits^{\scriptstyle}}\kern -1mm& H_{T}^{\bullet}
(X_{n,\overline{\eta}}^{T},{\Bbb Q}_{\ell})\cr}
$$
o\`{u} les fl\`{e}ches horizontales sont les fl\`{e}ches de
sp\'{e}cialisation pour $X_{n}$ et le $S$-sch\'{e}ma des points fixes
$X_{n}^{T}$ et o\`{u} les fl\`{e}ches verticales sont les fl\`{e}ches
de restriction au ferm\'{e} $X_{n}^{T}\subset X$. D'o\`{u} la
conclusion par passage \`{a} la limite inductive sur $n$.

\hfill\hfill$\square$

\section{4}{Un calcul d'homologie \'{e}quivariante}

On se propose de calculer l'homologie $\ell$-adique ${\Bbb G}_{{\rm
m},k}$-\'{e}quivariante de la cha\^{\i}ne de droites projectives \`{a}
la puissance $r$-\`{e}me $V^{r}$ pour l'action diagonale de ${\Bbb
G}_{{\rm m},k}\subset {\Bbb G}_{{\rm m},k}^{r}$.  On a d\'{e}j\`{a}
remarqu\'{e} que
$$
(V^{r})^{{\Bbb G}_{{\rm m},k}}= (V^{r})^{{\Bbb G}_{{\rm m},k}^{r}}
\cong {\Bbb Z}^{r} \times\mathop{\rm Spec}(k).
$$

\thm LEMME 4.1
\enonce
La fl\`{e}che de restriction
$$
{\Bbb Q}_{\ell}[{\Bbb Z}^{r}][t]= H_{\bullet}^{{\Bbb G}_{{\rm m},k}}
((V^{r})^{{\Bbb G}_{{\rm m},k}},{\Bbb Q}_{\ell})
\rightarrow H_{\bullet}^{{\Bbb G}_{{\rm m}, k}}(V^{r},
{\Bbb Q}_{\ell})
$$
est surjective et a pour noyau
$$
\bigoplus_{d=0}^{r-1}\Bigl(\bigcap_{{\scriptstyle R\subset
[1,r]\atop\scriptstyle |R|=r-d}}\sum_{\rho\in
R}(1-\tau_{\rho}){\Bbb Q}_{\ell}[{\Bbb Z}^{r}]\Bigr)t^{d}
$$
o\`{u} $[1,r]=\{1,\ldots ,r\}$, o\`{u} $|R|$ est le cardinal de $R$ et
o\`{u} $\tau_{\rho}$ est le $\rho$-\`{e}me \'{e}l\'{e}ment de la base
canonique de ${\Bbb Z}^{r}$, vu comme un \'{e}l\'{e}ment de ${\Bbb
Q}_{\ell}[{\Bbb Z}^{r}]$, de sorte que ${\Bbb Q}_{\ell}[{\Bbb
Z}^{r}]={\Bbb Q}_{\ell}[\tau_{1},\ldots
,\tau_{r},(\tau_{1}\cdots\tau_{r})^{-1}]$.

Par suite, on a
$$
H_{\bullet}^{{\Bbb G}_{{\rm m}, k}}(V^{r},{\Bbb
Q}_{\ell}) \cong \Bigl(\bigoplus_{d=0}^{r-1}{\Bbb Q}_{\ell}
(\Delta_{d})t^{d}\Bigr)\oplus {\Bbb Q}_{\ell}[{\Bbb Z}^{r}][t]t^{r}
$$
en tant que ${\Bbb Q}_{\ell}[{\Bbb Z}^{r}]$-module gradu\'{e}, o\`{u}
$$
\Delta_{d}=\bigcup_{{\scriptstyle R\subset [1,r]\atop\scriptstyle
|R|=r-d}}\{(\tau_{1},\ldots ,\tau_{r})\mid \tau_{\rho}=1,~\forall
\rho\in R\}\subset {\Bbb G}_{{\rm m},{\Bbb Q}_{\ell}}^{r}=\mathop{\rm
Spec}({\Bbb Q}_{\ell}[{\Bbb Z}^{r}])
$$
et o\`{u} ${\Bbb Q}_{\ell}(\Delta_{d})$ est le ${\Bbb Q}_{\ell}[{\Bbb
Z}^{r}]$-module des fonctions r\'{e}guli\`{e}res sur ${\Bbb
Q}_{\ell}$-sch\'{e}ma $\Delta_{d}$.
\endthm

\rem Preuve
\endrem
Commen\c{c}ons par calculer l'homologie $\ell$-adique ${\Bbb G}_{{\rm
m},k}$-\'{e}quivariante d'une copie $V$ de la cha\^{\i}ne
de droites projectives.  Comme $V$ est r\'{e}union croissante de
$k$-sch\'{e}mas projectifs \`{a} cohomologie $\ell$-adique pure, sur
lesquels ${\Bbb G}_{{\rm m}, k}$ agit avec un nombre fini de points
fixes et un nombre fini d'orbites de dimension $1$, la fl\`{e}che de
restriction
$$
{\Bbb Q}_{\ell}[{{\Bbb Z}}][t]=H_{\bullet}^{{\Bbb G}_{{\rm m},k}}
(\mathop{\rm Spec}(k)\times {\Bbb Z},{\Bbb Q}_{\ell})\rightarrow
H_{\bullet}^{{\Bbb G}_{{\rm m},k}}(V,{\Bbb Q}_{\ell}).
$$
est surjective et on peut calculer $H_{\bullet}^{{\Bbb G}_{{\rm m},k}}
(V,{\Bbb Q}_{\ell})$ par la m\'{e}thode de Goresky, Kottwitz et
MacPherson (calcul qu'ils ont d'ailleurs fait dans [1]). On trouve
que
$$
H_{\bullet}^{{\Bbb G}_{{\rm m},k}}(V,{\Bbb Q}_{\ell})={\Bbb
Q}_{\ell}[{{\Bbb Z}}][t]/(1-\tau ){\Bbb Q}_{\ell}[{{\Bbb Z}}]
$$
en tant que ${\Bbb Q}_{\ell}[{\Bbb Z}][\delta ]$-module, et donc que
$$
H_{\bullet}^{{\Bbb G}_{{\rm m},k}}(V,{\Bbb Q}_{\ell})\cong
{\Bbb Q}_{\ell}\oplus t{\Bbb Q}_{\ell}[{{\Bbb Z}}][t]
$$
en tant que ${\Bbb Q}_{\ell}[{\Bbb Z}]$-module, o\`{u}:
\vskip 1mm

\itemitem{-} $\tau$ est le g\'{e}n\'{e}rateur $1\in {\Bbb Z}$ vu comme
\'{e}l\'{e}ment de ${\Bbb Q}_{\ell}[{\Bbb Z}]$, de sorte que ${\Bbb
Q}_{\ell}[{\Bbb Z}]={\Bbb Q}_{\ell}[\tau ,\tau^{-1}]$,
\vskip 1mm

\itemitem{-} ${\Bbb Q}_{\ell}[{\Bbb Z}]\subset {\Bbb Q}_{\ell}[{{\Bbb
Z}}][t]$ est le noyau de $\delta$,
\vskip 1mm

\itemitem{-} dans la deuxi\`{e}me \'{e}galit\'{e}, le facteur direct
${\Bbb Q}_{\ell}$ est vu comme ${\Bbb Q}_{\ell}[{\Bbb Z}]$-Module via
l'homomorphisme ${\Bbb Q}_{\ell}[{\Bbb Z}]\rightarrow {\Bbb Q}_{\ell}$
qui envoie $\tau$ sur $1$.
\vskip 1mm

En passant \`{a} la puissance tensorielle $r$-\`{e}me (sur ${\Bbb
Q}_{\ell}$), on voit \`{a} l'aide de la formule de K\"{u}nneth que la
fl\`{e}che de restriction
$$
{\Bbb Q}_{\ell}[{\Bbb Z}^{r}][t_{1},\ldots ,t_{r}]=
H_{\bullet}^{{\Bbb G}_{{\rm m},k}^{r}}((V^{r})^{{\Bbb G}_{{\rm m},k}^{r}},
{\Bbb Q}_{\ell})\rightarrow H_{\bullet}^{{\Bbb G}_{{\rm m},k}^{r}}
(V^{r},{\Bbb Q}_{\ell})
$$
est aussi surjective et que
$$
H_{\bullet}^{{\Bbb G}_{{\rm m},k}^{r}}(V^{r},{\Bbb
Q}_{\ell})={\Bbb Q}_{\ell}[{\Bbb Z}^{r}][t_{1},\ldots ,t_{r}]/
\sum_{\rho =1}^{r} (1-\tau_{\rho}){\Bbb Q}_{\ell} [{\Bbb Z}^{r}]
[t_{1},\ldots ,\widehat{t}_{\rho},\ldots ,t_{r}].
$$

Comme la cohomologie $\ell$-adique de $V^{r}$ est pure,
la fl\`{e}che de restriction
$$
{\Bbb Q}_{\ell}[{\Bbb Z}^{r}][t]=H_{\bullet}^{{\Bbb G}_{{\rm m},k}}
((V^{r})^{{\Bbb G}_{{\rm m},k}},{\Bbb Q}_{\ell})\rightarrow
H_{\bullet}^{{\Bbb G}_{{\rm m},k}}(V^{r},{\Bbb Q}_{\ell})
$$
est surjective et n'est autre que la partie de la fl\`{e}che de
restriction en cohomologie ${\Bbb G}_{{\rm m}}^{r}$-\'{e}quivariante
annul\'{e}e par le noyau de l'\'{e}pimorphisme
$$
{\Bbb Q}_{\ell}[\delta_{1},\ldots ,\delta_{r}]\twoheadrightarrow {\Bbb
Q}_{\ell}[\delta ],~\delta_{\rho}\mapsto\delta ,~\forall \rho
=1,\ldots ,r,
$$
avec en particulier
$$\eqalign{
H_{\bullet}^{{\Bbb G}_{{\rm m},k}} ((V^{r})^{{\Bbb G}_{{\rm
m},k}},{\Bbb Q}_{\ell})={\Bbb Q}_{\ell}[{\Bbb Z}^{r}][t]
&\hookrightarrow {\Bbb Q}_{\ell}[{\Bbb Z}^{r}][t_{1},\ldots ,t_{r}]=
H_{\bullet}^{{\Bbb G}_{{\rm m},k}^{r}}((V^{r})^{{\Bbb G}_{{\rm
m},k}^{r}}, {\Bbb Q}_{\ell}).\cr
t&\mapsto t_{1}+\cdots +t_{r}\cr}
$$
Le noyau de la fl\`{e}che de restriction ${\Bbb Q}_{\ell}[{\Bbb
Z}^{r}][t]=H_{\bullet}^{{\Bbb G}_{{\rm m},k}} ((V^{r})^{{\Bbb G}_{{\rm
m},k}},{\Bbb Q}_{\ell}) \rightarrow H_{\bullet}^{{\Bbb G}_{{\rm
m},k}}(V^{r},{\Bbb Q}_{\ell})$ est donc form\'{e} des polyn\^{o}mes
$P(t)\in {\Bbb Q}_{\ell} [{\Bbb Z}^{r}][t]$ tels que
$$
P(t_{1}+\cdots +t_{r})\in\sum_{\rho =1}^{r}
(1-\tau_{\rho}){\Bbb Q}_{\ell} [{\Bbb Z}^{r}][t_{1},\ldots
,\widehat{t}_{\rho},\ldots ,t_{r}]
$$
c'est-\`{a}-dire des polyn\^{o}mes
$$
P(t)=\bigoplus_{d=0}^{r-1}p_{d}(\tau_{1},\ldots ,\tau_{r})t^{d}
$$
o\`{u}
$$
p_{d}(\tau_{1},\ldots ,\tau_{r})\in\bigcap_{{\scriptstyle R\subset
[1,r]\atop\scriptstyle |R|=r-d}}\sum_{\rho\in R}(1-\tau_{\rho}){\Bbb
Q}_{\ell}[\tau_{1},\ldots ,\tau_{r},(\tau_{1}\cdots\tau_{r})^{-1}].
$$
Par suite, on a
$$
H_{\bullet}^{{\Bbb G}_{{\rm m},k}}(V^{r},{\Bbb Q}_{\ell})={\Bbb
Q}_{\ell}[{\Bbb Z}^{r}][t]/\bigoplus_{d=0}^{r-1}\Bigl(
\bigcap_{{\scriptstyle R\subset [1,r]\atop\scriptstyle |R|=r-d}}
\sum_{\rho\in R}(1-\tau_{\rho}){\Bbb Q}_{\ell}[{\Bbb Z}^{r}]
\Bigr)t^{d}
$$
en tant que ${\Bbb Q}_{\ell}[{\Bbb Z}^{r}][\delta ]$-module, et donc
$$
H_{\bullet}^{{\Bbb G}_{{\rm m},k}}(V^{r},{\Bbb Q}_{\ell})\cong
\Bigl(\bigoplus_{d=0}^{r-1}{\Bbb Q}_{\ell}(\Delta_{d})t^{d}\Bigr) \oplus
{\Bbb Q}_{\ell}[{\Bbb Z}^{r}][t]t^{r}
$$
en tant que ${\Bbb Q}_{\ell}[{\Bbb Z}^{r}]$-module, d'o\`{u} le lemme.

\hfill\hfill$\square$
\vskip 3mm

\section{5}{Homologie \'{e}quivariante des fibres de Springer}

Pour calculer la cohomologie $H^{\ast}(Z_{I},{\cal L})$ qui intervient
dans le lemme fondamental g\'{e}om\'{e}trique, on suit la m\'{e}thode
propos\'{e}e par Goresky, Kottwitz et MacPherson, et on commence par
\'{e}tudier l'homologie \'{e}quivariante
$$
H_{\bullet}^{{\Bbb G}_{{\rm m},k}}(X_{I},{\Bbb Q}_{\ell})=
H_{\bullet}^{{\Bbb G}_{{\rm m},k}}(X_{I}^{0},{\Bbb Q}_{\ell})[{\Bbb Z}]
$$
pour l'action du sous-tore d\'{e}fini par la partition $I_{1}\amalg
I_{2}$,
$$
{\Bbb G}_{{\rm m},k}\cong {\Bbb G}_{{\rm m},k}^{2}/{\Bbb G}_{{\rm
m},k}\subset {\Bbb G}_{{\rm m},k}^{I}/{\Bbb G}_{{\rm
m},k}=T_{I},~u\mapsto (t_{i})_{i\in I},
$$
o\`{u} $t_{i}=u$ si $i\in I_{1}$ et $t_{i}=1$ si $i\in I_{2}$.

On consid\`{e}re pour cela le diagramme commutatif de ${\Bbb
Q}_{\ell}[X^{\ast}(\widetilde{T}_{\overline{\eta}})][\delta ]$-modules
gradu\'{e}s
$$\diagram{
H_{\bullet}^{{\Bbb G}_{{\rm m},\overline{\eta}}}
((\widetilde{\overline{P}}{}_{\overline{\eta}}^{\,\natural ,0})^{{\Bbb G}_{{\rm
m},\overline{\eta}}},{\Bbb Q}_{\ell}) &\kern -1mm\smash{\mathop{\hbox
to 8mm{\rightarrowfill}} \limits^{\scriptstyle}}\kern -1mm
&H_{\bullet}^{{\Bbb G}_{{\rm m},\overline{\eta}}}
((\overline{P}{}_{\overline{\eta}}^{\,\natural ,0})^{{\Bbb G}_{{\rm
m},\overline{\eta}}},{\Bbb Q}_{\ell}) &\kern -1mm\smash{\mathop{\hbox
to 8mm{\rightarrowfill}}\limits^{\scriptstyle}}\kern -1mm&
H_{\bullet}^{{\Bbb G}_{{\rm m},s}} ((\overline{P}{}_{s}^{\,\natural ,0})^{{\Bbb
G}_{{\rm m},s}},{\Bbb Q}_{\ell})\cr
\llap{$\scriptstyle $}\left\downarrow
\vbox to 4mm{}\right.\rlap{}&&\llap{$\scriptstyle $}\left\downarrow
\vbox to 4mm{}\right.\rlap{}&&\llap{}\left\downarrow
\vbox to 4mm{}\right.\rlap{$\scriptstyle $}\cr
H_{\bullet}^{{\Bbb G}_{{\rm m},\overline{\eta}}}
(\widetilde{\overline{P}}{}_{\overline{\eta}}^{\,\natural ,0}, {\Bbb Q}_{\ell})
&\kern -1mm\smash{\mathop{\hbox to 8mm {\rightarrowfill}}
\limits^{\scriptstyle }}\kern -1mm& H_{\bullet}^{{\Bbb G}_{{\rm m},
\overline{\eta}}} (\overline{P}{}_{\overline{\eta}}^{\,\natural ,0}, {\Bbb
Q}_{\ell})&\kern -1mm\smash{\mathop{\hbox to 8mm
{\rightarrowfill}}\limits^{\scriptstyle }}\kern -1mm&
H_{\bullet}^{{\Bbb G}_{{\rm m},s}} (\overline{P}{}_{s}^{\,\natural ,0},{\Bbb
Q}_{\ell})\cr}
$$
o\`{u} les fl\`{e}ches horizontales de gauche sont les fl\`{e}ches de
restriction par le morphisme
$\widetilde{\overline{P}}{}_{\overline{\eta}}^{\,\natural ,0}
\rightarrow \overline{P}{}_{\overline{\eta}}^{\,\natural ,0}$ et le
morphisme correspondant entre les ferm\'{e}s des points fixes sous
${\Bbb G}_{{\rm m},\overline{\eta}}$, o\`{u} les fl\`{e}ches
horizontales de droites sont les fl\`{e}ches de sp\'{e}cialisation
introduites dans la section $3$ pour le $S$-sch\'{e}ma
$\overline{P}{}^{\,\natural ,0}$, o\`{u} les fl\`{e}ches verticales
sont les fl\`{e}ches de restriction aux ferm\'{e}s des points fixes
sous ${\Bbb G}_{{\rm m}}$ et o\`{u} ${\Bbb
Q}_{\ell}[X^{\ast}(\widetilde{T}_{\overline{\eta}})]$ agit \`{a}
travers son quotient ${\Bbb Q}_{\ell}[X^{\ast}(T)]$ sur les
deuxi\`{e}me et troisi\`{e}me colonnes.

Compte tenu des lemmes 2.2 et 2.3, du fait que la cohomologie
$\ell$-adique est insensible aux hom\'{e}omorphismes et aux
changements de corps de base alg\'{e}briquement clos, et de la formule
de K\"{u}nneth, le carr\'{e} ext\'{e}rieur du diagramme ci-dessus se
r\'{e}crit
$$\diagram{
H_{\bullet}(X_{I_{1}}^{0}\times_{k}X_{I_{2}}^{0},{\Bbb Q}_{\ell})
[{\Bbb Z}^{r}][t] &\kern -14mm\smash{\mathop{\hbox to
21mm{\rightarrowfill}} \limits^{\scriptstyle}}\kern -1mm &H_{\bullet}
(X_{I_{1}}^{0} \times_{k}X_{I_{2}}^{0},{\Bbb Q}_{\ell})[{\Bbb
Z}][t]\cr
\llap{$\scriptstyle $}\left\downarrow
\vbox to 4mm{}\right.\rlap{}&&\llap{$\scriptstyle $}\left\downarrow
\vbox to 4mm{}\right.\rlap{}\cr
H_{\bullet}(X_{I_{1}}^{0}\times_{k}X_{I_{2}}^{0},{\Bbb Q}_{\ell})
\otimes_{{\Bbb Q}_{\ell}}H_{\bullet}^{{\Bbb G}_{{\rm m},k}}
(V^{r},{\Bbb Q}_{\ell}) &\kern -1mm
\smash{\mathop{\hbox to 16mm {\rightarrowfill}}\limits^{\scriptstyle }}
\kern -9mm& H_{\bullet}^{{\Bbb G}_{{\rm m},k}}(X_{I}^{0},{\Bbb Q}_{\ell})
\cr}\leqno{(\ast )}
$$
o\`{u} la fl\`{e}che verticale de gauche est le produit tensoriel
par $H_{\bullet}(X_{I_{1}}^{0}\times_{k}X_{I_{2}}^{0},{\Bbb
Q}_{\ell})$ de la fl\`{e}che de restriction
$$
{\Bbb Q}_{\ell}[{\Bbb Z}^{r}][t]=H_{\bullet}^{{\Bbb G}_{{\rm m},k}}
((V^{r})^{{\Bbb G}_{{\rm m},k}},{\Bbb Q}_{\ell})
\rightarrow H_{\bullet}^{{\Bbb G}_{{\rm m},k}}
(V^{r},{\Bbb Q}_{\ell})
$$
du lemme 4.1, o\`{u} la fl\`{e}che horizontale du haut est
l'\'{e}pimorphisme induit par le morphisme somme ${\Bbb Z}^{r}
\twoheadrightarrow {\Bbb Z}$ et o\`{u} la fl\`{e}che verticale de
droite est la fl\`{e}che de restriction
$$
H_{\bullet} (X_{I_{1}}^{0} \times_{k}X_{I_{2}}^{0},{\Bbb
Q}_{\ell})[{\Bbb Z}][t]=H_{\bullet}^{{\Bbb G}_{{\rm m},k}}
((X_{I}^{0})^{{\Bbb G}_{{\rm m},k}},{\Bbb Q}_{\ell})\rightarrow
H_{\bullet}^{{\Bbb G}_{{\rm m},k}}(X_{I}^{0},{\Bbb Q}_{\ell}).
$$

D'apr\`{e}s le lemme $4.1$, la fl\`{e}che verticale de gauche du
carr\'{e} $(\ast )$ est surjective et a pour noyau
$$\displaylines{
\qquad\widetilde{N}_{\overline{\eta},\bullet}=\bigoplus_{d=0}^{r-1}\Bigl(
\bigcap_{{\scriptstyle R\subset [1,r]\atop\scriptstyle |R|=r-d}}
\sum_{\rho\in R}(1-\tau_{\rho})H_{\bullet} (X_{I_{1}}^{0}
\times_{k}X_{I_{2}}^{0},{\Bbb Q}_{\ell}) [{\Bbb Z}^{r}]\Bigr)t^{d}
\hfill\cr\hfill
\subset H_{\bullet}(X_{I_{1}}^{0}\times_{k}X_{I_{2}}^{0},{\Bbb Q}_{\ell})
[{\Bbb Z}^{r}][t].\qquad}
$$

\thm LEMME 5.1
\enonce
L'image
$$
N_{\overline{\eta},\bullet}\subset H_{\bullet}(X_{I_{1}}^{0}
\times_{k}X_{I_{2}}^{0},{\Bbb Q}_{\ell}) [{\Bbb Z}][t]
$$
de $\widetilde{N}_{\overline{\eta},\bullet}$ par la fl\`{e}che
horizontale du haut du carr\'{e} $(\ast )$ est \'{e}gale \`{a}
$$
N_{\overline{\eta},\bullet}=\bigoplus_{d=0}^{r-1}(1-\tau )^{{r\choose
d}}H_{\bullet} (X_{I_{1}}^{0}\times_{k}X_{I_{2}}^{0},{\Bbb
Q}_{\ell})[{\Bbb Z}]t^{d}.
$$
\endthm

\rem Preuve
\endrem
L'\'{e}pimorphisme ${\Bbb Q}_{\ell}[{\Bbb Z}^{r}]\twoheadrightarrow
{\Bbb Q}_{\ell}[{\Bbb Z}]$ donn\'{e} par la somme ${\Bbb Z}^{r}
\rightarrow {\Bbb Z}$ n'est autre que
$$
{\Bbb Q}_{\ell}[\tau_{1},\ldots ,\tau_{r},(\tau_{1}\cdots\tau_{r})^{-1})
\twoheadrightarrow {\Bbb Q}_{\ell}[\tau ,\tau^{-1}],~\tau_{\rho}
\mapsto\tau ,~\forall \rho =1,\ldots ,r.
$$

On raisonne alors en terme de sch\'{e}mas affines sur ${\Bbb
Q}_{\ell}$ plut\^{o}t que de ${\Bbb Q}_{\ell}$-alg\`{e}bres.  Pour
chaque entier $d\geq 0$, le ferm\'{e} de ${\Bbb G}_{{\rm m},{\Bbb
Q}_{\ell}}^{r}$ d\'{e}fini par l'id\'{e}al
$$
\bigcap_{{\scriptstyle R\subset [1,r]\atop\scriptstyle |R|=r-d}}
\sum_{\rho\in R}(1-\tau_{\rho}){\Bbb Q}_{\ell}[{\Bbb Z}^{r}]
$$
est la r\'{e}union sur les $R\subset [1,r]$ tels que $|R|=r-d$ des
ferm\'{e}s
$$
\{(\tau_{1},\ldots ,\tau_{r})\mid \tau_{\rho}=1,~\forall\rho\in R\}.
$$
L'intersection de cette r\'{e}union de ferm\'{e}s avec la diagonale
$$
\tau_{\rho}=\tau ,~\forall\rho =1,\ldots ,r,
$$
est donc le point $\tau =1$ compt\'{e} avec multiplicit\'{e}
${r\choose r-d}={r\choose d}$.

\hfill\hfill$\square$
\vskip 3mm

\thm PROPOSITION 5.2
\enonce
Le noyau $N_{s,\bullet}$ de la fl\`{e}che de restriction
$$
H_{\bullet}(X_{I_{1}}^{0}\times_{k}X_{I_{2}}^{0},{\Bbb Q}_{\ell})
[{\Bbb Z}][t]= H_{\bullet}^{{\Bbb G}_{{\rm m},k}}((X_{I}^{0})^{{\Bbb
G}_{{\rm m},k}},{\Bbb Q}_{\ell}) \rightarrow H_{\bullet}^{{\Bbb
G}_{{\rm m},k}}(X_{I}^{0},{\Bbb Q}_{\ell})
$$
est born\'{e} inf\'{e}rieurement par
$$
\bigoplus_{d=0}^{r-1}(1-\tau )^{{r\choose d}}H_{\bullet}(X_{I_{1}}^{0}
\times_{k}X_{I_{2}}^{0},{\Bbb Q}_{\ell})[{\Bbb Z}]t^{d}\subset
N_{s,\bullet}\subset H_{\bullet}(X_{I_{1}}^{0}\times_{k}
X_{I_{2}}^{0},{\Bbb Q}_{\ell})[{\Bbb Z}][t].
$$
\endthm

\rem Preuve
\endrem
On a
$$
N_{\overline{\eta},\bullet}\subset N_{s,\bullet}\subset
H_{\bullet}(X_{I_{1}}^{0}\times_{k}X_{I_{2}}^{0},{\Bbb Q}_{\ell})
[{\Bbb Z}][t]
$$
puisque le carr\'{e} $(\ast )$ est commutatif.

\hfill\hfill$\square$
\vskip 3mm

{\it Nous supposerons dans la suite que la conjecture de puret\'{e}
$1.1$ est v\'{e}rifi\'{e}e}.  Alors, la fl\`{e}che de restriction en
homologie $\ell$-adique ${\Bbb G}_{{\rm m},k}$-\'{e}quivariante
$$
H_{\bullet}(X_{I_{1}}^{0}\times_{k}X_{I_{2}}^{0},{\Bbb Q}_{\ell})
[{\Bbb Z}][t]= H_{\bullet}^{{\Bbb G}_{{\rm m},k}}((X_{I}^{0})^{{\Bbb
G}_{{\rm m},k}},{\Bbb Q}_{\ell}) \rightarrow H_{\bullet}^{{\Bbb
G}_{{\rm m},k}}(X_{I}^{0},{\Bbb Q}_{\ell})
$$
est surjective.  Pour calculer $H_{\bullet}^{{\Bbb G}_{{\rm
m},k}}(X_{I}^{0},{\Bbb Q}_{\ell})$, il suffit donc de calculer
$N_{s,\bullet}$.

\section{6}{Le cas particulier $|I|=2$}

Supposons dans cette section que $I=\{i,j\}$ et que
$I_{1}=\{i\}$ et $I_{2}=\{j\}$.  Dans ce cas, $X_{I_{1}}^{0}=Z_{i}$ et
$X_{I_{2}}^{0}=Z_{j}$ sont des $k$-sch\'{e}mas projectifs.

\thm PROPOSITION 6.1
\enonce
Supposons que la conjecture de puret\'{e} $1.1$ est v\'{e}rifi\'{e}e.
Alors le noyau $N_{s,\bullet}$ est born\'{e} inf\'{e}rieurement et
sup\'{e}rieurement par
$$
\bigoplus_{d=0}^{r-1}(1-\tau )^{{r\choose d}}H_{\bullet}(Z_{i}
\times_{k}Z_{j},{\Bbb Q}_{\ell})[{\Bbb Z}]t^{d}\subset
N_{s,\bullet}\subset \bigoplus_{d=0}^{r-1}H_{\bullet}(Z_{i}\times_{k}
Z_{j},{\Bbb Q}_{\ell})[{\Bbb Z}]t^{d}
$$
dans $H_{\bullet}(Z_{i}\times_{k}Z_{j},{\Bbb Q}_{\ell})[{\Bbb Z}][t]$
\endthm

On remarquera que le quotient de la borne sup\'{e}rieure par la borne
inf\'{e}rieure est un ${\Bbb Q}_{\ell}[{\Bbb Z}]$-module de torsion
tu\'{e} par une puissance de $1-\tau$ qui ne d\'{e}pend que de $r$.

\rem Preuve
\endrem
Notons $D$ le rang total du ${\Bbb Q}_{\ell}$-espace vectoriel
$$
H_{\bullet}(Z_{i}\times_{k}Z_{j},{\Bbb Q}_{\ell}).
$$
Le ${\Bbb Q}_{\ell}[{\Bbb Z}]$-module
$$
\bigoplus_{d=0}^{r-1}(1-\tau )^{{r\choose d}}H_{\bullet}(Z_{i}
\times_{k}Z_{j},{\Bbb Q}_{\ell})[{\Bbb Z}]t^{d}
$$
est libre de rang $rD$ et on a
$$\eqalign{
\bigoplus_{d=0}^{r-1}(1-\tau )^{{r\choose d}}H_{\bullet}(Z_{i}
\times_{k}Z_{j},{\Bbb Q}_{\ell})[{\Bbb Z}]t^{d}\subset
N_{s,\bullet}\subset &\bigoplus_{d=0}^{r-1}H_{\bullet}(Z_{i}
\times_{k}Z_{j},{\Bbb Q}_{\ell})[{\Bbb Z}]t^{d}\cr
&\oplus
\bigoplus_{d=r}^{+\infty}H_{\bullet}(Z_{i} \times_{k}Z_{j},{\Bbb
Q}_{\ell})[{\Bbb Z}]t^{d}.\cr}
$$
Il suffit donc de d\'{e}montrer que le ${\Bbb Q}_{\ell}[{\Bbb
Z}]$-module $N_{s,\bullet}$ est lui aussi libre de rang $rD$.
Comme sous-${\Bbb Q}_{\ell}[{\Bbb Z}]$-module d'un module libre,
$N_{s,\bullet}$ est sans torsion et il suffit m\^{e}me de
d\'{e}montrer qu'il est de type fini et de rang $rD$.

Pour tout $k$-sch\'{e}ma de type fini $a:S\rightarrow \mathop{\rm
Spec}(k)$, il est commode d'introduire la cohomologie \'{e}quivariante
$$
R\Gamma_{{\Bbb G}_{{\rm m}},{\rm c}}^{\bullet} (S,{\Bbb Q}_{\ell})=
R\Gamma_{{\Bbb G}_{{\rm m}}}^{\bullet} (\mathop{\rm Spec}(k),
R\overline{a}_{!}{\Bbb Q}_{\ell})
$$
o\`{u} $\overline{a}:[S/{\Bbb G}_{{\rm m}}]\rightarrow B({\Bbb
G}_{{\rm m}})$ est le morphisme de champs induit par $a$, et
l'homologie \'{e}quivariante
$$
H_{\bullet}^{{\Bbb G}_{{\rm m}},{\rm c}} (S,{\Bbb Q}_{\ell})=
\mathop{\rm Hom}\nolimits_{{\Bbb Q}_{\ell}}(H_{{\Bbb G}_{{\rm m}},{\rm
c}}^{\bullet}(S,{\Bbb Q}_{\ell}),{\Bbb Q}_{\ell}).
$$
On se gardera de confondre ces groupes de (co)homologie avec ceux de
la (co)homologie \`{a} supports compacts du $k$-champ alg\'{e}brique
$[S/{\Bbb G}_{{\rm m}}]$; on esp\`{e}re que la notation adopt\'{e}e
ici ne pr\^{e}te pas \`{a} confusion.  Ces d\'{e}finitions
s'\'{e}tendent aux $k$-sch\'{e}mas localement de type fini qui
interviennent ici.  Avec cette notion d'homologie, on a
$$
N_{\bullet}=H_{\bullet +1}^{{\Bbb G}_{{\rm m}},{\rm c}}
(X_{I}^{0}-(X_{I}^{0})^{{\Bbb G}_{{\rm m}}},{\Bbb Q}_{\ell}).
$$

Pour all\'{e}ger les notations dans la suite de la d\'{e}monstration,
nous supprimerons l'indice $s$ qui ne joue plus aucun r\^{o}le, et
aussi l'indice $k$ du produit fibr\'{e} $\times_{k}$ et la
r\'{e}f\'{e}rence \`{a} ${\Bbb Q}_{\ell}$ dans les groupes de
(co)homologie.

Consid\'{e}rons alors le fibr\'{e} vectoriel de rang $r$
$$
U\rightarrow \coprod_{\lambda\in {\Bbb Z}}
(X_{i}^{\lambda}\times X_{j}^{-\lambda})\cong {\Bbb
Z}\times (Z_{i}\times Z_{j})
$$
dont la fibre en $(M_{i}\subset E_{i},M_{j}\subset E_{j})$ est
l'espace vectoriel
$$
\mathop{\rm Hom}\nolimits_{{\cal O}_{F}[x]}(M_{j}, E_{i}/M_{i})=
\mathop{\rm Hom}\nolimits_{{\cal O}_{F}[x]}(M_{j}/P_{i}
(\gamma_{j})M_{j}, P_{j}(\gamma_{i})^{-1}M_{i}/M_{i})
$$
o\`{u} $x$ agit sur $E_{i}$ et $M_{i}$ par multiplication par
$\gamma_{i}$ et agit sur $E_{j}$ et $M_{j}$ par multiplication par
$\gamma_{j}$, et o\`{u} les $k$-espaces vectoriels
$P_{j}(\gamma_{i})^{-1}M_{i}/M_{i}$ et $M_{j}/P_{i}(\gamma_{j})M_{j}$
sont tous deux de rang $r$ (cf.  [4]).

On a un morphisme bijectif, mais non radiciel,
$$
U\rightarrow X_{I}^{0}
$$
qui envoie $(M_{i},M_{j},\varphi\in \mathop{\rm Hom}\nolimits_{{\cal
O}_{F}[x]}(M_{j}, E_{i}/M_{i}))$ sur le ${\cal O}_{F}$-r\'{e}seau
$M\subset E_{I}$ dont la trace sur le facteur direct $E_{i}$ de
$E_{I}$ est $M_{i}$ et la projection sur le facteur direct $E_{j}$ de
$E_{I}$ est $M_{j}$, et qui tel que
$$
M/M_{i}\subset (E_{i}/M_{i})\oplus M_{j}
$$
soit le graphe de $\varphi$.  L'application non alg\'{e}brique
$X_{I}^{0}\rightarrow {\Bbb Z}\times (Z_{i}\times Z_{j})$
consid\'{e}r\'{e}e par Kazhdan et Lusztig dans la section 5 de [5]
n'est autre que la compos\'{e}e de l'inverse de cette bijection et de
la projection $U\rightarrow {\Bbb Z}\times (Z_{i}\times Z_{j})$.

Par construction, le fibr\'{e} vectoriel et le morphisme bijectif
ci-dessus sont $\Lambda_{I}^{0}$-\'{e}quivariants. De plus, on a
$$
\Bigl(U\rightarrow \coprod_{\lambda\in {\Bbb Z}}
(X_{i}^{\lambda}\times X_{j}^{-\lambda})\Bigr)\cong {\Bbb
Z}\times (U'\rightarrow Z_{i}\times Z_{j})
$$
o\`{u} $U'$ est la restriction de $U$ \`{a} $X_{i}^{0}\times
X_{j}^{0}$.

On stratifie $U'$ et $U$, par le rang de $\varphi\in \mathop{\rm
Hom}\nolimits_{{\cal O}_{F}[x]} (M_{j}/P_{i} (\gamma_{j})M_{j},
P_{j}(\gamma_{i})^{-1}M_{i}/M_{i})$.  On obtient une stratification
$U'=\bigcup_{\rho =0}^{r}U_{\rho}'$ en parties quasi-projectives sur
$k$ et localement ferm\'{e}es dans $U'$, et la stratification
$U=\bigcup_{\rho =0}^{r}({\Bbb Z}\times U_{\rho}')$

On consid\`{e}re alors la partition $X_{I}^{0}=\bigcup_{\rho
=0}^{r}S_{\rho}$ image de cette stratification de $U$ par le morphisme
bijectif $U\rightarrow X_{I}^{0}$ ci-dessus.  On a $S_{\rho}={\Bbb
Z}\times S_{\rho}'$ o\`{u} $S_{\rho}'$ est l'image de $U_{\rho}'$.

On v\'{e}rifie que les $S_{\rho}={\Bbb Z}\times S_{\rho}'$ sont des
parties localement ferm\'{e}es de $X_{I}^{0}$ et forment une
stratification de $X_{I}^{0}$.  On v\'{e}rifie aussi que les
stratifications de $U$ et $X_{I}^{0}$ que l'on vient d'introduire sont
$\Lambda_{I}^{0}$-\'{e}quivariantes.  On v\'{e}rifie en outre que
$S_{\rho}'$ est {\it isomorphe \`{a} $U_{\rho}'$ en tant que
$k$-sch\'{e}ma}.

Le groupe multiplicatif ${\Bbb G}_{{\rm m}}$ agit par homoth\'{e}ties
sur le fibr\'{e} vectoriel $U$.  On v\'{e}rifie que le morphisme
bijectif $U\rightarrow X_{I}^{0}$ et les stratifications de $U$ et
$X_{I}^{0}$ ci-dessus sont ${\Bbb G}_{{\rm m}}$-\'{e}quivariants, et
que les strates $U_{0}$ et $S_{0}$ sont \'{e}gales aux ferm\'{e}s
$U^{{\Bbb G}_{{\rm m}}}$ et $X^{{\Bbb G}_{{\rm m}}}$ et sont donc
toutes les deux isomorphes \`{a} ${\Bbb Z}\times (Z_{i}\times Z_{j})$.

Des stratifications $X_{I}^{0}-(X_{I}^{0})^{{\Bbb G}_{{\rm m}}}=
\bigcup_{\rho =1}^{r}({\Bbb Z}\times S_{\rho}')$ et $U'-U'^{{\Bbb
G}_{{\rm m}}}= \bigcup_{\rho =1}^{r}U_{\rho}'$, on d\'{e}duit des
suites spectrales en homologie $\ell$-adique ${\Bbb G}_{{\rm
m}}$-\'{e}quivariante
$$
E_{\rho,\sigma}^{1}=H_{\rho +\sigma}^{{\Bbb G}_{{\rm m}},{\rm c}}
(S_{\rho +1}')[{\Bbb Z}]\Rightarrow H_{\rho +\sigma}^{{\Bbb G}_{{\rm
m}},{\rm c}}(X_{I}^{0}-(X_{I}^{0})^{{\Bbb G}_{{\rm m}}})
$$
et
$$
F_{\rho,\sigma}^{1}=H_{\rho +\sigma}^{{\Bbb G}_{{\rm m}},{\rm c}}
(U_{\rho +1}')[{\Bbb Z}]\Rightarrow H_{\rho +\sigma}^{{\Bbb G}_{{\rm
m}},{\rm c}}(U'-U'^{{\Bbb G}_{{\rm m}}})[{\Bbb Z}]
$$
avec
$$
E_{\rho,\sigma}^{1}=F_{\rho,\sigma}^{1}.
$$

Ces suites spectrales sont des suites spectrales de ${\Bbb
Q}_{\ell}[{\Bbb Z}]$-modules dont les termes communs $E^{1}=F^{1}$
sont libres de type fini.  On a donc montr\'{e} que le ${\Bbb
Q}_{\ell} [{\Bbb Z}]$-module $N_{\bullet}=H_{\bullet +1}^{{\Bbb
G}_{{\rm m}},{\rm c}} (X_{I}^{0}-(X_{I}^{0})^{{\Bbb G}_{{\rm m}}})$
est de type fini et que, pour chaque entier $w$, on a
$$
\sum_{n}(-1)^{n}\mathop{\rm rang}([N_{n}]_{w})= \sum_{n}(-1)^{n}
\mathop{\rm dim}([H_{n}^{{\Bbb G}_{{\rm m}},{\rm c}}(U'-U'^{{\Bbb
G}_{{\rm m}}})]_{w}),
$$
o\`{u} $\mathop{\rm rang}(\cdot )$ d\'{e}signe le rang
g\'{e}n\'{e}rique d'un ${\Bbb Q}_{\ell}[{\Bbb Z}]$-module de type
fini, alors que $\mathop{\rm dim}(\cdot )$ d\'{e}signe la dimension
d'un ${\Bbb Q}_{\ell}$-espace vectoriel, et o\`{u} $[\cdot ]_{w}$
d\'{e}signe la partie de poids $w$ d'un groupe de (co)homologie
$\ell$-adique.

Or, d'apr\`{e}s le lemme $6.2$ ci-dessous, on a
$$
H_{{\Bbb G}_{{\rm m}},{\rm c}}^{\bullet +1} (U'-U'^{{\Bbb G}_{{\rm
m}}})=H^{\bullet}(Z_{i}\times Z_{j})[\delta ]/
(\delta^{r}-c_{1}\delta^{r-1}+\cdots +(-1)^{r}c_{r})
$$
pour des classes $c_{d}\in H^{2d}(Z_{i}\times Z_{j})(d)$,
et donc par division euclidienne, on a un isomorphisme canonique de
${\Bbb Q}_{\ell}$-espaces vectoriels
$$
H_{{\Bbb G}_{{\rm m}},{\rm c}}^{n+1} (U'-U'^{{\Bbb G}_{{\rm
m}}}) \cong \bigoplus_{d=0}^{r-1}H^{n-2d}(Z_{i}\times Z_{j})\delta^{d}
$$
pour chaque entier $n$.

Finallement, la conjecture $1.1$ assure que, pour chaque entier $n$,
$$
N_{n}\subset\bigoplus_{d=0}^{r-1}H_{n-2d}(Z_{i}\times Z_{j})[{\Bbb Z}]
t^{d}
$$
et $H_{n+1}^{{\Bbb G}_{{\rm m}},{\rm c}}(U'-U'^{{\Bbb G}_{{\rm m}}})$
sont purs de poids $n$.

On a donc montr\'{e} que
$$
\mathop{\rm rang}(N_{n})=\mathop{\rm dim}(H_{n+1}^{{\Bbb
G}_{{\rm m}},{\rm c}}(U'-U'^{{\Bbb G}_{{\rm m}}}))=
\sum_{d=0}^{r-1}\mathop{\rm dim}(H_{n-2d}(Z_{i}\times Z_{j}))
$$
pour chaque entier $n$, et donc que
$$
\sum_{n}\mathop{\rm rang}(N_{n})=r\sum_{n} \mathop{\rm
dim}(H_{n-2d}(Z_{i}\times Z_{j}))=rD,
$$
ce qui termine la preuve de la proposition.
\hfill\hfill$\square$
\vskip 3mm

\thm LEMME 6.2
\enonce
Soit $S$ un $k$-sch\'{e}ma connexe projectif et $E={\Bbb V}({\cal
E})\rightarrow S$ un fibr\'{e} vectoriel de rang $r$.  Soit
$E^{\circ}\subset E$ l'ouvert compl\'{e}mentaire de la section nulle,
muni de son action naturelle de ${\Bbb G}_{{\rm m},k}$.  Alors
$$
H_{{\Bbb G}_{{\rm m},k},{\rm c}}^{\bullet +1}(E^{\circ},{\Bbb Q}_{\ell})=
H^{\bullet}(S,{\Bbb Q}_{\ell})[\delta ]/(\delta^{r}-c_{1}({\cal
E})\delta^{r-1}+\cdots +(-1)^{r}c_{r}({\cal E}))
$$
en tant que $H^{\bullet}(S,{\Bbb Q}_{\ell})[\delta ]$-module gradu\'{e}.
\endthm

\rem Preuve
\endrem
On a
$$
R\Gamma_{{\Bbb G}_{{\rm m},k},{\rm c}}(E^{\circ},{\Bbb Q}_{\ell})=
R\Gamma_{{\Bbb G}_{{\rm m},k}}(S,R\overline{p}_{\ast}
R\overline{q}_{!}{\Bbb Q}_{\ell})
$$
o\`{u} on a factoris\'{e} la projection canonique $E^{\circ}
\rightarrow S$ en $E^{\circ}\,\smash{\mathop{\hbox to 6mm
{\rightarrowfill}} \limits^{\scriptstyle q}}\,{\Bbb P}(E) \,
\smash{\mathop{\hbox to 6mm{\rightarrowfill}} \limits^{\scriptstyle
p}}\,S$ et o\`{u} on a not\'{e} $\overline{q}:{\Bbb P} ({\cal E})
=[E^{\circ}/{\Bbb G}_{{\rm m},k}]\rightarrow [{\Bbb P}({\cal E})/
{\Bbb G}_{{\rm m},k}]$ et $\overline{p}:[{\Bbb P}({\cal E})/{\Bbb
G}_{{\rm m},k}] \rightarrow [S/{\Bbb G}_{{\rm m},k}]$ les morphismes
correspondants entre champs quotients.  Or on a un triangle
distingu\'{e}
$$
{\Bbb Q}_{\ell,[{\Bbb P}({\cal E})/{\Bbb G}_{{\rm m},k}]}[-2](-1)\rightarrow
{\Bbb Q}_{\ell,[{\Bbb P}({\cal E})/{\Bbb G}_{{\rm m},k}]}\rightarrow
R\overline{q}_{!}{\Bbb Q}_{\ell}[1]\rightarrow
$$
o\`{u} la premi\`{e}re fl\`{e}che est le cup-produit par $c_{1}({\cal
O}_{{\Bbb P}({\cal E})}(1))-\delta$.  En effet, on v\'{e}rifie que la
premi\`{e}re classe de Chern du fibr\'{e} en droites sur le champ
quotient $[{\Bbb P}({\cal E})/{\Bbb G}_{{\rm m},k}]$ dont le
compl\'{e}mentaire de la section nulle est $\overline{q}:{\Bbb P}
({\cal E}) =[E^{\circ}/{\Bbb G}_{{\rm m},k}]\rightarrow [{\Bbb P}
({\cal E})/ {\Bbb G}_{{\rm m},k}]$ est pr\'{e}cis\'{e}ment
$c_{1}({\cal O}_{{\Bbb P}({\cal E})}(1))-\delta$ dans $H_{{\Bbb
G}_{{\rm m},k}}^{2}({\Bbb P}({\cal E}),{\Bbb Q}_{\ell})(1)=H^{2}({\Bbb
P}({\cal E}),{\Bbb Q}_{\ell})(1)\oplus H^{0}({\Bbb P}({\cal E}),{\Bbb
Q}_{\ell})\delta$.

Or, on sait que l'on a
$$
R\overline{p}_{\ast}{\Bbb Q}_{\ell}=\bigoplus_{d=0}^{r-1}
R\overline{p}_{\ast}^{2d}{\Bbb Q}_{\ell}[-2d]={\Bbb Q}_{\ell ,S}
[\Delta ]/(\Delta^{r}-c_{1}({\cal E})
\Delta^{r-1}+\cdots +(-1)^{r}c_{r}({\cal E}))
$$
o\`{u} $\Delta$ est en fait $c_{1}({\cal O}_{{\Bbb P}({\cal E})}(1))$.
Par suite on a un triangle distingu\'{e}
$$\displaylines{
\qquad R\Gamma (S,{\Bbb Q}_{\ell})[\Delta ,\delta ]/
(\Delta^{r}-c_{1}({\cal E})\Delta^{r-1}+\cdots +(-1)^{r}c_{r}({\cal E}))
\hfill\cr\hfill
\rightarrow R\Gamma (S,{\Bbb Q}_{\ell})[\Delta ,\delta ]/
(\Delta^{r}-c_{1}({\cal E})\Delta^{r-1}+\cdots +(-1)^{r}c_{r}({\cal E}))
\hfill\cr\hfill
\rightarrow R\Gamma_{{\Bbb G}_{{\rm m},k},{\rm c}}(E^{\circ},{\Bbb
Q}_{\ell})[1]\rightarrow\qquad}
$$
o\`{u} la premi\`{e}re fl\`{e}che est la multiplication par $\Delta
-\delta$, d'o\`{u} la formule annonc\'{e}e
$$
R\Gamma_{{\Bbb G}_{{\rm m},k},{\rm c}}(E^{\circ},{\Bbb
Q}_{\ell})[1]=R\Gamma (S,{\Bbb Q}_{\ell})[\delta ]/
(\delta^{r}-c_{1}({\cal E})\delta^{r-1}+\cdots +(-1)^{r}c_{r}({\cal E})).
$$
\hfill\hfill$\square$

\section{7}{Le cas g\'{e}n\'{e}ral}

Nous allons utiliser la m\^{e}me m\'{e}thode de r\'{e}duction au rang
$1$ que dans [1].

\thm LEMME 7.1
\enonce
Pour chaque partie $J$ de $I$, le lieu des points fixes
$$
X_{J}^{T_{J}}=\prod_{i\in J}X_{i}
$$
est le $k$-sch\'{e}ma des r\'{e}seaux $M\subset E_{J}$ de la forme
$M=\bigoplus_{i\in J}M_{i}$ pour des r\'{e}seaux $M_{i}\subset E_{i}$
tels que $\gamma_{i}M_{i}\subset M_{i}$.

La r\'{e}union de ce lieu des points fixes et des $T_{J}$-orbites de
dimension $1$ dans $X_{J}$ est le ferm\'{e}
$$
X_{J,1}\subset X_{J}
$$
r\'{e}union sur les sous-ensembles $\{i,j\}\subset J$ \`{a} deux
\'{e}l\'{e}ments, des ferm\'{e}s
$$
X_{J}^{ij}=X_{J}^{\mathop{\rm Ker}(\chi_{ij})}
=X_{ij}\times\Bigl(\prod_{k\in J-\{i,j\}}X_{k}\Bigr)\subset X_{J}
$$
o\`{u} $\chi_{ij}:T_{J}\rightarrow {\Bbb G}_{{\rm m},k}$ est le
caract\`{e}re qui envoie $(t_{k})_{k\in I}$ sur $t_{i}/t_{j}$,
ferm\'{e}s form\'{e}s des r\'{e}seaux $M$ qui se d\'{e}composent en
$M=M_{ij}\oplus\bigoplus_{{\scriptstyle k\in J\atop\scriptstyle
k\not=i,j}}M_{k}$ pour des r\'{e}seaux $M_{ij}\subset E_{ij}$ et
$M_{k}\subset E_{k}$ tels que $\gamma_{ij}M_{ij}\subset M_{ij}$
et $\gamma_{k}M_{k}\subset M_{k}$.

De plus, on a
$$
X_{J}^{ij}\cap X_{J}^{i'j'}=\prod_{i\in J}X_{i}
$$
si $\{i,j\}\not=\{i',j'\}$.

\hfill\hfill$\square$
\endthm

Notons ${\frak m}\subset {\Bbb Q}_{\ell}[{\Bbb Z}^{I}]$ l'id\'{e}al
maximal noyau de l'homomorphisme de ${\Bbb Q}_{\ell}$-alg\`{e}bre qui
envoie $\tau_{i}$ sur $1$ si $i\in I_{1}$ et sur $-1$ si $i\in I_{2}$.

\thm PROPOSITION 7.2
\enonce
Supposons que la conjecture de puret\'{e} $1.1$ est v\'{e}rifi\'{e}e.
Alors on a l'inclusion
$$
H_{\bullet+1}^{T_{I},{\rm c}}(X_{I}-(X_{I_{1}}\times_{k}
X_{I_{2}},{\Bbb Q}_{\ell}))_{{\frak m}}\subset H_{\bullet}^{T_{I}}
(X_{I_{1}}\times_{k} X_{I_{2}},{\Bbb Q}_{\ell})_{{\frak m}}
\left\{\prod_{i\in I_{1},j\in I_{2}}
(\delta_{i}-\delta_{j})^{r_{ij}}\right\}
$$
entre sous-${\Bbb Q}_{\ell}[{\Bbb Z}^{I}]_{{\frak
m}}[(\delta_{i})_{i\in I}]$-modules de $H_{\bullet}^{T_{I}}
(X_{I_{1}}\times_{k}X_{I_{2}},{\Bbb Q}_{\ell})_{{\frak m}}$.
\endthm

Dans cet \'{e}nonc\'{e}, $(\cdot )_{{\frak m}}=(\cdot )\otimes_{{\Bbb
Q}_{\ell}[{\Bbb Z}^{I}]}{\Bbb Q}_{\ell}[{\Bbb Z}^{I}]_{{\frak m}}$ est
la localisation en l'id\'{e}al maximal ${\frak m}$, et on a
utilis\'{e} la notation de [1]
$$
M\{P\}=\{m\in M\mid Pm=(0)\}
$$
pour tout ${\Bbb Q}_{\ell}[(\delta_{i})_{i\in I}]$-Module $M$ et tout
\'{e}l\'{e}ment $P\in {\Bbb Q}_{\ell}[(\delta_{i})_{i\in I}]$.

\rem Preuve
\endrem
Nous all\'{e}gerons de nouveau les notations en supprimant l'indice
$k$ des produits fibr\'{e}s $\times_{k}$ et la r\'{e}f\'{e}rence \`{a}
${\Bbb Q}_{\ell}$ dans les groupes de cohomologie.

D'apr\`{e}s conjecture de puret\'{e} $1.1$, on a le diagramme commutatif
$$\diagram{
0&&0\cr
\llap{$\scriptstyle $}\left\downarrow
\vbox to 4mm{}\right.\rlap{}&&\llap{}\left\downarrow
\vbox to 4mm{}\right.\rlap{$\scriptstyle $}\cr
H_{\bullet +1}^{T_{I},{\rm c}}(X_{I_{1}}\times X_{I_{2}}
-X_{I}^{T_{I}})&\kern -1mm\smash{\mathop{\hbox to
8mm{\rightarrowfill}} \limits^{\scriptstyle }}\kern -1mm&H_{\bullet
+1}^{T_{I}, {\rm c}}(X_{I}-X_{I}^{T_{I}})\cr
\llap{$\scriptstyle $}\left\downarrow
\vbox to 4mm{}\right.\rlap{}&&\llap{}\left\downarrow
\vbox to 4mm{}\right.\rlap{$\scriptstyle $}\cr
H_{\bullet}^{T_{I}}(X_{I}^{T_{I}})&\kern
-12mm\vbox{\hrule width 24mm\vskip2pt\hrule width 24mm}\kern
-6mm&H_{\bullet}^{T_{I}}(X_{I}^{T_{I}})\cr
\llap{$\scriptstyle $}\left\downarrow
\vbox to 4mm{}\right.\rlap{}&&\llap{}\left\downarrow
\vbox to 4mm{}\right.\rlap{$\scriptstyle $}\cr
H_{\bullet}^{T_{I}}(X_{I_{1}}\times X_{I_{2}})
&\kern -7mm\smash{\mathop{\hbox to 20mm{\rightarrowfill}}
\limits_{\scriptstyle }}\kern -7mm&H_{\bullet}^{T_{I}}(X_{I})\cr
\llap{$\scriptstyle $}\left\downarrow
\vbox to 4mm{}\right.\rlap{}&&\llap{}\left\downarrow
\vbox to 4mm{}\right.\rlap{$\scriptstyle $}\cr
0&&0\cr}
$$
o\`{u} les fl\`{e}ches horizontales sont les fl\`{e}ches de
restriction et o\`{u} toutes les homologies consid\'{e}r\'{e}es sont
pures de poids $n$ en chaque degr\'{e} $n$ d'apr\`{e}s la conjecture
$1.1$.  On a aussi une suite exacte longue
$$\eqalign{
&H_{\bullet +1}^{T_{I},{\rm c}}(X_{I_{1}}\times X_{I_{2}}
-X_{I}^{T_{I}})\rightarrow H_{\bullet +1}^{T_{I},{\rm c}}
(X_{I}-X_{I}^{T_{I}}) \rightarrow H_{\bullet +1}^{T_{I},{\rm c}}
(X_{I}-(X_{I_{1}}\times X_{I_{2}}))\cr
\rightarrow &H_{\bullet}^{T_{I},{\rm
c}}(X_{I_{1}}\times X_{I_{2}} -X_{I}^{T_{I}}),\cr}
$$
dont la premi\`{e}re fl\`{e}che est injective et a une source et un
but qui sont purs de poids $n$ en chaque degr\'{e} $n$.  Cette suite
exacte longue se r\'{e}duit donc \`{a} une suite exacte courte
$$
0\rightarrow H_{\bullet +1}^{T_{I},{\rm c}}(X_{I_{1}}\times X_{I_{2}}
-X_{I}^{T_{I}})\rightarrow H_{\bullet +1}^{T_{I},{\rm c}}
(X_{I}-X_{I}^{T_{I}}) \rightarrow H_{\bullet +1}^{T_{I},{\rm c}}
(X_{I}-(X_{I_{1}}\times X_{I_{2}}))\rightarrow 0
$$
o\`{u} tous les termes sont purs de poids $n$ en chaque degr\'{e} $n$,
et le lemme du serpent nous donne alors une suite exacte qui n'est
autre que la suite exacte courte
$$
0\rightarrow H_{\bullet +1}^{T_{I},{\rm c}}(X_{I}-(X_{I_{1}}\times
X_{I_{2}}))\rightarrow H_{\bullet}^{T_{I}}(X_{I_{1}}
\times X_{I_{2}}) \rightarrow H_{\bullet}^{T_{I}}(X_{I})\rightarrow
0
$$
de restriction au ferm\'{e} $X_{I_{1}}\times X_{I_{2}}$ de $X_{I}$.

Toujours d'apr\`{e}s la conjecture de puret\'{e}, pour chaque partie
$J$ de $I$ on a en outre la suite exacte (4.3.1) de [1]
$$
H_{\bullet +1}^{T_{J},{\rm c}}(X_{J,1}-X_{J}^{T_{J}}) \rightarrow
H_{\bullet}^{T_{J}}(X_{J}^{T_{J}}) \rightarrow
H_{\bullet}^{T_{J}}(X_{J}) \rightarrow 0\,;
$$
en d'autres termes, la fl\`{e}che de restriction
$$
H_{\bullet +1}^{T_{J},{\rm c}}(X_{J,1}-X_{J}^{T_{J}})
\rightarrow H_{\bullet +1}^{T_{J},{\rm c}}(X_{J}-X_{J}^{T_{J}})
$$
est surjective.

On a donc un morphisme de suite exactes courtes
$$
\matrix{H_{\bullet +1}^{T_{I},{\rm c}}(X_{I_{1},I_{2},1}-
X_{I}^{T_{I}}) & \hookrightarrow & H_{\bullet +1}^{T_{I},{\rm c}}
(X_{I,1}-X_{I}^{T_{I}}) & \twoheadrightarrow & H_{\bullet
+1}^{T_{I},{\rm c}} (X_{I,1}-X_{I_{1},I_{2},1}) \cr
\downarrow &  & \downarrow &  & \downarrow \cr
H_{\bullet +1}^{T_{I},{\rm c}}(X_{I_{1}}\times X_{I_{2}}
-X_{I}^{T_{I}}) & \hookrightarrow & H_{\bullet +1}^{T_{I},{\rm c}}
(X_{I}-X_{I}^{T_{I}}) & \twoheadrightarrow & H_{\bullet
+1}^{T_{I},{\rm c}}(X_{I}-X_{I_{1}}\times X_{I_{2}})\cr}
$$
o\`{u} $X_{I_{1},I_{2},1}=X_{I,1}\cap (X_{I_{1}}\times
X_{I_{2}})=(X_{I_{1},1}\times X_{I_{2}}^{T_{I_{2}}})\cup
(X_{I_{1}}^{T_{I_{1}}}\times X_{I_{2},1})$ et o\`{u} les deux
premi\`{e}res fl\`{e}ches verticales, et par suite aussi la
troisi\`{e}me, sont surjectives.

Comme
$$
H_{\bullet}^{T_{J}}(X_{J}^{T_{J}})=H_{\bullet}^{T_{J}}\Bigl(\prod_{i\in
J}X_{i}^{0}\Bigr)[{\Bbb Z}^{J}]=H_{\bullet}(X_{J}^{T_{J}})[(t_{i})_{i\in J}]
=H_{\bullet}\Bigl(\prod_{i\in J}X_{i}^{0}\Bigr)[{\Bbb Z}^{J}] 
[(t_{i})_{i\in J}]
$$
et
$$\eqalign{
H_{\bullet}^{T_{J},{\rm c}}(X_{J,1}-X_{J}^{T_{J}})&=
\bigoplus_{{\scriptstyle \{i,j\}\subset J\atop\scriptstyle
i\not=j}}H_{\bullet}^{T_{J},{\rm c}}(X_{J}^{ij}-X_{J}^{T_{J}})\cr
&=\bigoplus_{{\scriptstyle \{i,j\}\subset J\atop\scriptstyle
i\not=j}}H_{\bullet}^{T_{ij},{\rm c}}(X_{ij}-X_{ij}^{T_{ij}})\cr
&\kern 15mm\otimes H_{\bullet}\Bigl(\prod_{k\in
J-\{i,j\}}X_{k}^{0}\Bigr)[{\Bbb Z}^{J-\{i,j\}}][(t_{k})_{k\in J-\{i,j\}}],\cr}
$$
on a une surjection
$$\displaylines{
\qquad\bigoplus_{{\scriptstyle i\in I_{1}\atop\scriptstyle
j\in I_{2}}}H_{\bullet +1}^{T_{ij},{\rm c}}(X_{ij}-X_{ij}^{T_{ij}})
\otimes H_{\bullet}\Bigl(\prod_{k\in J-\{i,j\}}X_{k}^{0}\Bigl)[{\Bbb
Z}^{I-\{i,j\}}][(t_{k})_{k\in I-\{i,j\}}]
\hfill\cr\hfill
\twoheadrightarrow H_{\bullet +1}^{T_{I},{\rm c}}
(X_{I}-X_{I_{1}}\times X_{I_{2}}).\qquad}
$$

Soient $i\in I_{1}$ et $j\in I_{2}$.  En utilisant la proposition
$6.1$ pour $I:=\{i,j\}$, on obtient les inclusions
$$\displaylines{
\qquad\bigoplus_{d=0}^{r_{ij}-1}(\tau_{i}-\tau_{j})^{{r_{ij}\choose
d}} H_{\bullet}(X_{i}^{0}\times X_{j}^{0})[{\Bbb Z}^{\{i,j\}}]
[t_{i}+t_{j}] (t_{i}-t_{j})^{d}
\hfill\cr\hfill
\subset H_{\bullet +1}^{T_{ij},{\rm c}}(X_{ij}^{0}-(X_{i,j}^{0})^{T_{ij}})
\subset\bigoplus_{d=0}^{r_{ij}-1}H_{\bullet}(X_{i}^{0}\times X_{j}^{0})
[{\Bbb Z}^{\{i,j\}}][t_{i}+t_{j}](t_{i}-t_{j})^{d}\qquad}
$$
de sous-${\Bbb Q}_{\ell}[{\Bbb Z}^{\{i,j\}}]$-modules de
$$
H_{\bullet}(X_{i}^{0}\times X_{j}^{0})[{\Bbb Z}^{\{i,j\}}]
[t_{i},t_{j}]=H_{\bullet}(X_{i}^{0}\times X_{j}^{0})[{\Bbb
Z}^{\{i,j\}}][t_{i}+t_{j},t_{i}-t_{j}].
$$

Comme l'\'{e}l\'{e}ment $\tau_{i}-\tau_{j}$ est inversible
dans l'anneau local ${\Bbb Q}_{\ell}[{\Bbb Z}^{I}]_{{\frak m}}$, les
inclusions ci-dessus donnent, apr\`{e}s localisation en ${\frak m}$,
l'\'{e}galit\'{e}
$$\displaylines{
\qquad \left(H_{\bullet +1}^{T_{ij},{\rm c}}(X_{ij}-X_{i,j}^{T_{ij}})
\otimes H_{\bullet}\Bigl(\prod_{k\in I-\{i,j\}}X_{k}^{0}\Bigr)[{\Bbb
Z}^{I-\{i,j\}}][(t_{k})_{k\in I-\{i,j\}}]\right)_{{\frak m}}
\hfill\cr\hfill
\eqalign{&=\bigoplus_{d=0}^{r_{ij}-1}H_{\bullet}(X_{I}^{T_{I}})_{{\rm m}}
[(t_{k})_{k\in I-\{i,j\}}][t_{i}+t_{j}](t_{i}-t_{j})^{d}\cr
&=H_{\bullet}(X_{I}^{T_{I}})_{{\rm m}}[(t_{k})_{k\in I}]
\{(\delta_{i}-\delta_{j})^{r_{ij}}\}.\cr}\qquad}
$$
\vskip 3mm

En r\'{e}sum\'{e}, les diagrammes commutatifs d'inclusions
$$
\matrix{(X_{ij}-X_{ij}^{T_{ij}})\times\prod_{k\in I-\{i,j\}}X_{k}
& \subset & X_{ij}\times\prod_{k\in I-\{i,j\}}X_{k} & \supset &
X_{ij}^{T_{ij}}\times\prod_{k\in I-\{i,j\}}X_{k}\cr
\cap &  & \cap &  & \cap\cr
X_{I,1}-X_{I_{1},I_{2},1} & \subset & X_{I,1} & \supset &
X_{I_{1},I_{2},1}\cr
\cap &  & \cap &  & \cap\cr
X_{I}-X_{I_{1}}\times X_{I_{2}} & \subset & X_{I} & \supset &
X_{I_{1}}\times X_{I_{2}}\cr}
$$
pour $i\in I_{1}$ et $j\in I_{2}$ induisent un carr\'{e} commutatif de
${\Bbb Q}_{\ell}[{\Bbb Z}^{I}]_{{\frak m}}[(\delta_{i})_{i\in
I}]$-modules
$$\diagram{
\bigoplus_{{\scriptstyle i\in I_{1}\atop\scriptstyle
j\in I_{2}}}H_{\bullet}(X_{I}^{T_{I}})_{{\rm m}}[(t_{k})_{k\in I}]
\{(\delta_{i}-\delta_{j})^{r_{ij}}\}
&\kern -1mm\smash{\mathop{\hbox to 8mm{\rightarrowfill}}
\limits^{\scriptstyle }}\kern -1mm&H_{\bullet}(X_{I}^{T_{I}})_{{\frak m}}
[(t_{k})_{k\in I}]\cr
\llap{$\scriptstyle $}\left\downarrow
\vbox to 4mm{}\right.\rlap{}&&\llap{}\left\downarrow
\vbox to 4mm{}\right.\rlap{$\scriptstyle $}\cr
H_{\bullet+1}^{T_{I},{\rm c}}(X_{I}-(X_{I_{1}}\times X_{I_{2}}))_{{\frak m}}
&\kern -10mm\smash{\mathop{\hbox to 17mm{\rightarrowfill}}
\limits_{\scriptstyle }}\kern -1mm&H_{\bullet}^{T_{I}}
(X_{I_{1}}\times X_{I_{2}})_{{\frak m}}\cr}
$$
o\`{u} les fl\`{e}ches verticales sont surjectives et o\`{u} la
fl\`{e}che horizontale du bas est injective. Or l'image
$$
\sum_{{\scriptstyle i\in I_{1}\atop\scriptstyle j\in I_{2}}}
H_{\bullet}(X_{I}^{T_{I}})_{{\rm m}}[(t_{k})_{k\in I}]
\{(\delta_{i}-\delta_{j})^{r_{ij}}\}\subset
H_{\bullet}(X_{I}^{T_{I}})_{{\frak m}} [(t_{k})_{k\in I}]
$$
de la fl\`{e}che horizontale du haut du carr\'{e} ci-dessus n'est
autre que
$$
H_{\bullet}(X_{I}^{T_{I}})_{{\rm m}}[(t_{k})_{k\in I}]
\textstyle\left\{\prod_{{\scriptstyle i\in I_{1}\atop\scriptstyle
j\in I_{2}}}(\delta_{i}-\delta_{j})^{r_{ij}}\right\}
$$
d'apr\`{e}s le lemme 3.2 de [1].  On a donc un carr\'{e} commutatif de
${\Bbb Q}_{\ell}[{\Bbb Z}^{I}]_{{\frak m}}[(\delta_{i})_{i\in
I}]$-modules
$$\diagram{
H_{\bullet}(X_{I}^{T_{I}})_{{\rm m}}[(t_{k})_{k\in I}]
\left\{\prod_{{\scriptstyle i\in I_{1}\atop\scriptstyle
j\in I_{2}}}(\delta_{i}-\delta_{j})^{r_{ij}}\right\}
&\kern -1mm\smash{\mathop{\lhook\joinrel\mathrel
{\hbox to 8mm{\rightarrowfill}}}\limits^{\scriptstyle }}
\kern -1mm&H_{\bullet}(X_{I}^{T_{I}})_{{\frak m}}[(t_{k})_{k\in I}]\cr
\llap{$\scriptstyle $}\left\downarrow
\vbox to 4mm{}\right.\rlap{}&&\llap{}\left\downarrow
\vbox to 4mm{}\right.\rlap{$\scriptstyle $}\cr
H_{\bullet+1}^{T_{I},{\rm c}}(X_{I}-(X_{I_{1}}\times
X_{I_{2}}))_{{\frak m}} &\kern -10mm\smash{\mathop{\lhook
\joinrel\mathrel {\hbox to17mm{\rightarrowfill}}}
\limits_{\scriptstyle }}\kern -1mm&H_{\bullet}^{T_{I}}
(X_{I_{1}}\times X_{I_{2}})_{{\frak m}}\cr}
$$
o\`{u} les fl\`{e}ches horizontales sont injectives et o\`{u} les
fl\`{e}ches verticales sont surjectives, d'o\`{u} la conclusion.
\hfill\hfill$\square$
\vskip 3mm

Passons maintenant \`{a} la cohomologie \'{e}quivariante pour le
sous-tore ${\Bbb G}_{{\rm m},k}^{2}$ de $T_{I}$ d\'{e}fini par
l'inclusion
$$
{\Bbb G}_{{\rm m},k}^{2}\hookrightarrow T_{I},~(u_{1},u_{2})\mapsto
(t_{i})_{i\in I}
$$
o\`{u} $t_{i}=u_{\alpha}$ si $i\in I_{\alpha}$, pour $\alpha =1,2$.
On a un \'{e}pimorphisme canonique
$$
{\Bbb Q}_{\ell}[(\delta_{i})_{i\in I}]=H_{T_{I}}^{\bullet}(\mathop{\rm
Spec}(k),{\Bbb Q}_{\ell})\twoheadrightarrow H_{{\Bbb G}_{{\rm
m},k}^{2}}^{\bullet}(\mathop{\rm Spec}(k),{\Bbb Q}_{\ell})=
{\Bbb Q}_{\ell}[\partial_{1},\partial_{2}]
$$
qui envoie $\delta_{i}$ sur $\partial_{\alpha}$ pour tout $i\in
I_{\alpha}$ et $\alpha =1,2$, \'{e}pimorphisme dont le noyau est
l'id\'{e}al ${\frak a}$ engendr\'{e} par $\delta_{i}-\delta_{j}$ pour
$i,j\in I_{\alpha}$ et $\alpha =1,2$.  D'apr\`{e}s Goresky, Kottwitz
et MacPherson, on a
$$
H_{\bullet}^{{\Bbb G}_{{\rm m},k}^{2}}(X_{I_{1}}\times_{k}X_{I_{2}},
{\Bbb Q}_{\ell})=H_{\bullet}^{T_{I}} (X_{I_{1}}\times X_{I_{2}},
{\Bbb Q}_{\ell})\{{\frak a}\}
$$
et
$$
H_{\bullet}^{{\Bbb G}_{{\rm m},k}^{2}}(X_{I},{\Bbb Q}_{\ell})=
H_{\bullet}^{T_{I}}(X_{I},{\Bbb Q}_{\ell})\{{\frak a}\}
$$
puisque les cohomologies $\ell$-adiques ordinaires de
$X_{I_{1}}\times_{k}X_{I_{2}}$ et $X_{I}$ sont pures.  Par suite, on
a aussi
$$
H_{\bullet+1}^{{\Bbb G}_{{\rm m},k}^{2},{\rm c}}(X_{I}-(X_{I_{1}}
\times_{k}X_{I_{2}}),{\Bbb Q}_{\ell})=H_{\bullet+1}^{T_{I},{\rm c}}
(X_{I}-(X_{I_{1}}\times_{k}X_{I_{2}}),{\Bbb Q}_{\ell})_{{\frak m}}
\{{\frak a}\}
$$
puisque l'on sait a priori que la suite
$$
0\rightarrow H_{\bullet+1}^{{\Bbb G}_{{\rm m},k}^{2},{\rm
c}}(X_{I}-(X_{I_{1}}\times_{k}X_{I_{2}}),{\Bbb Q}_{\ell})
\rightarrow H_{\bullet}^{{\Bbb G}_{{\rm m},k}^{2}}(X_{I},
{\Bbb Q}_{\ell})\rightarrow H_{\bullet}^{{\Bbb G}_{{\rm
m},k}^{2}}(X_{I_{1}}\times X_{I_{2}},{\Bbb Q}_{\ell})\rightarrow 0
$$
est exacte (on a une suite exacte longue \'{e}vidente et on utilise la
conjecture de puret\'{e} pour voir que la fl\`{e}che
$H_{\bullet}^{{\Bbb G}_{{\rm m},k}^{2}}(X_{I}, {\Bbb
Q}_{\ell})\rightarrow H_{\bullet}^{{\Bbb G}_{{\rm
m},k}^{2}}(X_{I_{1}}\times X_{I_{2}},{\Bbb Q}_{\ell})$ est
surjective).

Comme la diagonale de ${\Bbb G}_{{\rm m},k}^{2}$ agit trivialement, on
peut m\^{e}me remplacer ${\Bbb G}_{{\rm m},k}^{2}$ par son quotient ${\Bbb
G}_{{\rm m},k}^{2}/{\Bbb G}_{{\rm m},k}\cong {\Bbb G}_{{\rm m},k}$
dans ce qui pr\'{e}c\`{e}de sans autres effets que de remplacer
$(u_{1},u_{2})$ par $t=u_{1}-u_{2}$ et $(\partial_{1},\partial_{2})$
par $\delta =\partial_{1}-\partial_{2}$.

On d\'{e}duit donc de la proposition $7.2$ le corollaire:

\thm COROLLAIRE 7.3
\enonce
Supposons que la conjecture de puret\'{e} $1.1$ est v\'{e}rifi\'{e}e.
Alors, on a l'inclusion
$$
H_{\bullet+1}^{{\Bbb G}_{{\rm m},k},{\rm c}}(X_{I}-(X_{I_{1}}
\times_{k} X_{I_{2}}),{\Bbb Q}_{\ell})_{{\frak m}}\subset
H_{\bullet}(X_{I_{1}}\times_{k}X_{I_{2}},
{\Bbb Q}_{\ell})_{{\frak m}}[t]\{\delta^{r}\}
$$
entre sous-${\Bbb Q}_{\ell}[{\Bbb Z}^{I}]_{{\frak m}}[\delta
]$-modules de $H_{\bullet}^{{\Bbb G}_{{\rm m},k}} (X_{I_{1}}\times_{k}
X_{I_{2}},{\Bbb Q}_{\ell})_{{\frak m}}=H_{\bullet}(X_{I_{1}}
\times_{k}X_{I_{2}},{\Bbb Q}_{\ell})_{{\frak m}}[t]$, o\`{u} on
rappelle que $r=\sum_{i\in I_{1},j\in J_{2}}r_{ij}$.
\endthm

Mettant ensemble les r\'{e}sultats du corollaire 7.3 et de la
proposition 5.2, on obtient finalement:

\thm TH\'{E}OR\`{E}ME 7.4
\enonce
On a en fait l'\'{e}galit\'{e}
$$
H_{\bullet+1}^{{\Bbb G}_{{\rm m},k},{\rm c}}(X_{I}-(X_{I_{1}}
\times_{k} X_{I_{2}}),{\Bbb Q}_{\ell})_{{\frak m}}=
H_{\bullet}(X_{I_{1}}\times_{k}X_{I_{2}},
{\Bbb Q}_{\ell})_{{\frak m}}[t]\{\delta^{r}\}
$$
entre sous-${\Bbb Q}_{\ell}[{\Bbb Z}^{I}]_{{\frak m}}[\delta
]$-modules de $H_{\bullet}(X_{I_{1}} \times_{k}X_{I_{2}},{\Bbb
Q}_{\ell})_{{\frak m}}[t]$, de sorte que
$$
H_{\bullet}^{{\Bbb G}_{{\rm m},k}}(X_{I},{\Bbb Q}_{\ell})_{{\frak
m}}=t^{r}H_{\bullet}(X_{I_{1}} \times_{k}X_{I_{2}},{\Bbb
Q}_{\ell})_{{\frak m}}[t].
$$
\hfill\hfill$\square$
\endthm

On passe finalement au quotient par $\Lambda_{I}$ comme dans [1]
pour obtenir le th\'{e}or\`{e}me $1.3$.
\vskip 10mm

\newtoks\ref \newtoks\auteur \newtoks\titre \newtoks\annee
\newtoks\revue \newtoks\tome \newtoks\pages \newtoks\reste

\def\bibitem#1{\parindent=10pt\itemitem{#1}\parindent=24pt}

\def\article{\bibitem{[\the\ref]}%
\the\auteur~-- \the\titre, {\sl\the\revue} {\bf\the\tome},
({\the\annee}), \the\pages.\smallskip\filbreak}

\def\autre{\bibitem{[\the\ref]}%
\the\auteur~-- \the\reste.\smallskip\filbreak}

\ref={1}
\auteur={M. {\pc GORESKY}, R. {\pc KOTTWITZ}, R. {\pc MACPHERSON}}
\reste={Homology of affine Springer fibers in the unramified case,
http://www.math.ias.edu/$\sim$goresky/preprints.html, (2002)}
\autre

\ref={2}
\auteur={G. {\pc LAUMON}}
\reste={Fibres de Springer et jacobiennes compactifi\'{e}es,
http://arxiv.org/abs/ math.AG/0204109, (2002)}
\autre

\ref={3}
\auteur={G. {\pc LUSZTIG}, J.M. {\pc SMELT}}
\titre={Fixed point varieties on the space of lattices}
\revue={Bull. London. Math. Soc.}
\tome={23}
\annee={1991}
\pages={213-218}
\article

\ref={4}
\auteur={G. {\pc LAUMON}, M. {\pc RAPOPORT}}
\reste={A geometric approach to the fundamental lemma for unitary
groups, http://arxiv.org/abs/math.AG/9711021, (1997)}
\autre

\ref={5}
\auteur={D {\pc KAZHDAN}, G. {\pc LUSZTIG}}
\titre={Fixed Point Varieties on Affine Flag Manifolds}
\revue={Israel J. of Math.}
\tome={62}
\annee={1988}
\pages={129-168}
\article

\bye